\documentclass[sort&compress, preprint,12pt]{elsarticle}

\journal{}
\usepackage[flushleft]{threeparttable}
\usepackage{graphicx}
\graphicspath{{Figures/}}
\usepackage{graphics}
\usepackage{epsfig}
\usepackage{amsfonts}
\usepackage{amssymb}
\usepackage{amsmath}
\usepackage{commath}
\usepackage{amsopn}
\usepackage{hyperref}
\usepackage{cleveref}
\usepackage{tikz}
\usepackage{pgfplots}
\usepgfplotslibrary{groupplots}
\pgfplotsset{compat=newest}
\usepackage{filecontents}
\usepackage{lineno}
\usepackage{xspace}
\usepackage{array}
\usepackage[font = normal, labelfont = bf]{caption}
\usepackage[figurename = Fig.]{caption}
\usepackage{caption}
\usepackage{cite}
\usepackage{multirow}
\usepackage{enumerate}
\usepackage{mathrsfs}
\usepackage{pdflscape}
\usepackage[FIGBOTCAP]{subfigure}
\usepackage[utf8]{inputenc}
\usepackage{rotating}
\usepackage{empheq}
\usepackage{standalone}
\usepackage{listings}
\usepackage{geometry}
\usepackage{booktabs, caption, makecell}

\usepackage{threeparttable}
\usepackage{dcolumn}
\newcolumntype{d}[1]{D{.}{.}{#1}}
\usepackage{listings}
\usepackage{color}
\usepackage{natbib}

\begin{document}

\begin{frontmatter}

\title{An extended polygonal finite element method for large deformation fracture analysis}

\author[ath1]{Hai D. Huynh}
\author[ath2]{Phuong Tran}
\author[ath3]{Xiaoying Zhuang}
\author[ath1]{H. Nguyen-Xuan\corref{cor1}}

\address[ath1]{CIRTech Institute, Chi Minh City University of Technology (HUTECH), Ho Chi Minh City, Vietnam}
\address[ath2]{Department of Civil $\&$ Infrastructure Engineering, RMIT University, Melbourne VIC 3001, Australia}
\address[ath3]{Institute of Continuum Mechanics, Leibniz Universität Hannover, Appelstraße 11, 30167 Hannover, Germany}
%\cortext[cor1]{Corresponding author. Tel: (+084) 906682393.\\
%	           E-mail address: \itshape{ngx.hung@hutech.edu.vn}  (H. Nguyen-Xuan).}
\cortext[cor1]{Corresponding author. Tel: (+084) 906682393.}
\ead{ngx.hung@hutech.edu.vn}
%\fntext[ath1]{Email: \itshape{huynhdonghai13@gmail.com}}

\begin{abstract}
The modeling of large deformation fracture mechanics has been a challenging problem regarding the accuracy of numerical methods and their ability to deal with considerable changes in deformations of meshes where having the presence of cracks. This paper further investigates the extended finite element method (XFEM) for the simulation of large strain fracture for hyper-elastic materials, in particular rubber ones. A crucial idea is to use a polygonal mesh to represent space of the present numerical technique in advance, and then a local refinement of structured meshes at the vicinity of the discontinuities is additionally established. Due to differences in the size and type of elements at the boundaries of those two regions, hanging nodes produced in the modified mesh are considered as normal nodes in an arbitrarily polygonal element. Conforming these special elements becomes straightforward by the flexible use of basis functions over polygonal elements. Results of this study are shown through several numerical examples to prove its efficiency and accuracy through comparison with former achievements.

\end{abstract}


\end{frontmatter}

\section{Introduction}
\label{sec 1}

The extended finite element method (XFEM) has become a well-known tool to numerically simulate discontinuity problems. This technique was first proposed in \citep{Belytschko_Black_1999} to model issues related to cracks involved in strong discontinuities. Mo\"es and Sukumar \citep{Moes_Dolbow_1999,Daux_Moes_2000} then extended its applications to model structures with voids and close interfaces represented as weak discontinuities. The main feature of XFEM allows us generate the mesh background independently to the discontinuities by adding enriched parts to the classical finite approximations under the concept of the Partition of Unity \citep{Melenk_Babuska_1996}. The use of enrichment functions was accordingly taken into account with meshless methods \citep{Pant_Singh_Mishra_2010}, isogeometric analysis (IGA) \citep{Ghorash_Valizadeh_Mohammad_2011i}. In regard to the application in fracture mechanics, XFEM has been successfully employed in modeling crack propagation \citep{Duflot_2008,Fries_Belytschko_2010,Bayesteh_Mohammadi_2011}, cohesive fracture \citep{Mariani_Perego_2003,Unger_Eckardt_Konke}, cracks with multiple branches \citep{Daux_Moes_2000,Richardson_Hegemann}, three dimensional crack growth \citep{Baydoun_Fries_2012,Pedro_Belytschko_2005}, etc. With the two-dimensional analysis of XFEM, a domain is divided into non-overlapping subdomains, which are normally comprised of triangular or quadrilateral elements. Mesh discretization is still a challenging task because element types and element geometries significantly affect the approximation quality. In particular, while triangle elements suit to discretize domains with complicated shapes, they provide results with less accuracy than that of quadrilateral elements. In recent years, much attention is paid to applications of polygonal elements in numerical methods due to their ability to produce highly accurate solutions as well as their flexibility in mesh generation. The construction of polygonal finite elements was introduced in the work of Sukumar et al \citep{Sukumar_Tabarraei_2004,Sukumar_Malsch_2006}. Paulino et al \citep{Talischi_Paulino_Pereira_2012} then established a framework for two dimensional mesh generator written in Matlab. Their applications have been widely utilized in various fields, such as simulation of incompressible fluid flows \citep{Talischi_Pereira_Paulino_2014}, contact problems \citep{Biabanaki_2014,Khoei_2015}, three dimensional linear elasticity \citep{Gain_Talischi_Paulino_2014}, polycrystalline materials \citep{Ghosh_Moorthy_1995}, and hyper-elastic materials \citep{Chi_Talischi_Pamies_Paulino_2014}, scaled boundary polygonal finite element method (SBPFEM) \citep{Ooi_2017}, linear strain smoothing \citep{Pramod_2018}, advanced virtual element techniques (VEM) \citep{Chi_2017,Chi_2018}, limit analysis \citep{Hung_2017,Hung_2018} and so on.\\

The XFEM is extended from the standard FEM by the introduction of an enrichment technique for elements to capture the challenging features of fracture problems. The effects of a crack such as displacement jump and stress singularity near the crack tip are installed mathematically with proper enrichment functions \citep{Moes_Dolbow_1999} to incorporate into the displacement approximation. This approach avoids re-meshing during analysis of problems involved in discontinuities. However, XFEM still exists weaknesses regarding features of a typical FEM such as lower-order approximations and sensitivities to element distortion (in the case of large deformations). The development of enrichment approach has been applied to other numerical methods to overcome the mentioned difficulties. For example, the extended isogeometric analysis \citep{Ghorash_Valizadeh_Mohammad_2011i,Nhon_Valizadeh_2015,Ghorashi_Valizadeh_2015} characterizes high order approximations and high continuity but loses Kronecker-delta property, and extended meshfree methods \citep{Pant_Singh_Mishra_2010,Zhuang_Augarde_Bordas_2011} offer approximations based on space discretization by only nodes and not by elements. Early study on applying polygonal meshes into XFEM was presented in \citep{Tabarraei_Sukumar_2008,Khoei_Yasbo_2015} to model issues of linear elastic fracture. This topic is promising due to advantages and widespread applications of polygonal elements.\\

Most achievement in fracture mechanics is numerical investigations of linear failure analysis, in particular, for brittle materials being represented under a small deformation, while behaviors of ductile materials in this field are crucially necessary for research. Philippe \citep{Philippe_Wolfgang_1994_1,Philippe_Wolfgang_1994_2} examined finite strains at a crack tip for hyper-elastic Neo-Hookean materials. Several investigations for nonlinear failure analysis of constitutive laws in relation to cohesive models were developed by Dolbow \citep{Dolbow_Moes_Belytschko_2001}, Wells \citep{Wells_Sluysi_2001}, Rabczuk \citep{Rabczuk_Zi_2008}. They recently extended this work to simulate geometrically nonlinear fracture problems for nearly incompressible material \citep{Dolbow_Devan_2004}. A further approach to focus on the enrichment around crack tips was reported \citep{Legrain_Moes_Verron_2005} by Legrain to deal with large deformation problems in the context of XFEM. Besides, studies in the nonlinear fractures for compressible hyper-elastic materials have been achieved by Karoui \citep{Karoui_2014} investigated with the Ciarlet-Geymonat hyper-elastic model, and Rashetnia \citep{Rashetnia_2015} examined with classes of Neo-Hookean materials to simulate the finite deformation of cracks. A recent study on large deformation of ductile fracture problems \citep{Broumand_Khoei_2013} was investigated in the association with the non-local damage plasticity model. Such success has made contributions to the scientific community of the fracture mechanics field. However, these authors have examined the use of only common finite elements such as triangular or quadrilateral meshes in terms of the space discretization. The interest of polygonal finite element technique in modeling crack problems has recently introduced by Sukumar \citep{Tabarraei_Sukumar_2008}, which only focused on the linear fracture aspects. Thus, the extension of its advantages to the finite deformation in fracture mechanics will be definitely favorable.\\

The present work focuses on employment of XFEM with the polygonal element discretization to solve two-dimensional fracture mechanics problems with large deformation. Herein, the Neo-Hookean hyper-elastic models for both compressible and incompressible materials are used in the modified XFEM analysis. In the framework of polygonal finite elements, a correction gradient technique of shape functions is applied to enhance the accuracy of evaluating numerical integration of the weak form.  Due to the presence of a jump in displacements at vicinity of cracks, a structured mesh of quadrilateral elements is embedded in polygonal elements are embedded. The advantage from this approach is the creation of mesh refinements adapting into regions containing discontinuities. Moreover, the replacements are easy to control the mesh refinement around crack paths that may definitely bring improvements in the computational efficiency. It is observed that differences from the size as well as the type of elements along the boundary of two mesh regions lead to appearance of hanging nodes. However, such difficulty is completely solved through the flexibility of conforming interpolation functions employed in polygons. For the implementation, an updated Lagrangian formulation is taken into the modified XFEM. This approach allows the geometry of the computational domain to be updated periodically, and therefore tracking the deformation of cracks becomes efficiently. Meanwhile, the obstacle to transfer data of certain variables in the updated Lagrangian approach is facilitated by an additional mapping between the deformed element and its undeformed counterpart, found in \citep{Broumand_Khoei_2013}.\\%In addition, the data transfer of certain variables from the former mesh to the latter one is avoided as well.\\

The paper is outlined as follows: the formulations of non-linear fracture in XFEM based on the updated Lagrangian description are described in Section \ref{sec 2}. In Section \ref{sec 3}, the polygonal finite elements are employed in the present numerical method, followed by the construction of basis functions and the procedure of treating hanging nodes. Finally, several numerical examples are given in Section \ref{sec 4} to demonstrate the accuracy and effectiveness of the present work.

\section{Governing equations}
\label{sec 2}
\subsection{Large deformation formulation in FEM}
\label{sub_sec 2.1}

\begin{figure}[!httb]
	\centering
	\includegraphics[scale=1]{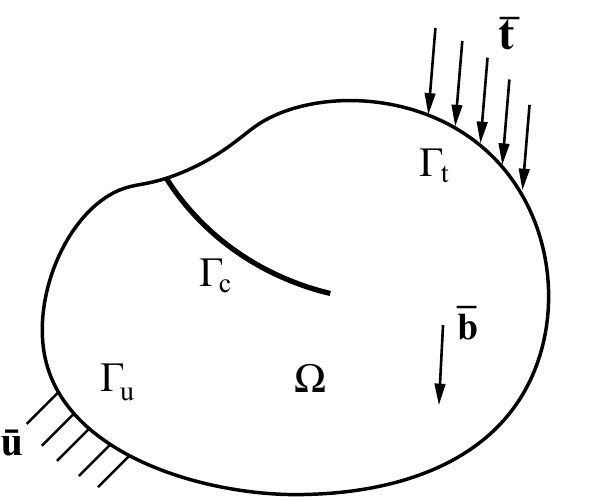}
	\caption{A crack problem with boundary conditions.}
	\label{F1}
\end{figure}

We consider a two-dimension body $\Omega$ bounded by $\Gamma$, which is subjected to a body force $\mathbf{\bar {b}}$. The boundary conditions are given by the prescribed traction $\mathbf{\bar {t}}$, and displacement $\mathbf{\bar{u}}$ on ${\Gamma _t}$ and ${\Gamma _u}$, respectively, and the crack surface ${\Gamma _c}$ is denoted by as illustrated in Fig. \ref{F1}

The equilibrium equation and boundary conditions defined in the updated Lagrangian description are given by

\begin{subequations}
	\begin{align}
	\nabla\cdot\boldsymbol{\sigma}  + \mathbf{\bar {b}} = \mathbf{0}  	&\quad \text{ in } 	\Omega	\\
	\boldsymbol{\sigma}\cdot\mathbf{n}= \mathbf{\bar{t}} 	&\quad	\text{ on } 	\Gamma _{t}	\\
	\boldsymbol{\sigma}\cdot\mathbf{n} = \mathbf{0}	&\quad	 \text{ on } 	\Gamma _c	\\
	\mathbf{u} = \mathbf{\bar{u}}	&\quad	 \text{ on } 	\Gamma _u
	\end{align}
	\label{Eq1}
\end{subequations}
where $\boldsymbol{\sigma}$ is the Cauchy stress tensor and $\mathbf{n}$ is the unit outward normal vector.\\

Following the finite element procedure, the weak form of the equilibrium equations is constructed by multiplying the test function $\delta {\bf{u}}$ and applying partial integration over the domain $\varphi (\Omega )$ which is transformed to the current configuration
\begin{equation}
\int\limits_{\varphi (\Omega )} {{\boldsymbol{\sigma}}:{\nabla _x}\delta {\bf{u}}} dv - \int\limits_{\varphi (\Omega )} {\delta {\bf{u}} \cdot {\bf{\bar b}}} dv - \int\limits_{\varphi ({\Gamma _t})} {\delta {\bf{u}} \cdot {\bf{\bar t}}} d\Gamma  = 0
\label{Eq4}
\end{equation}

In the first term of Eq. \ref{Eq4}, the Cauchy stress tensor with respect to the current configuration facilitates the gradient of the test function $\delta {\bf{u}}$ by its symmetric part with the following definition

\begin{equation}
\nabla _x^S\delta {\bf{u}} = \frac{1}{2}\left( {{\nabla _x}\delta {\bf{u}} + \nabla _x^T\delta {\bf{u}}} \right)
\label{Eq5}
\end{equation}
Hence, the weak form is rewritten as
\begin{equation}
\int\limits_{\varphi (\Omega )} {{\boldsymbol{\sigma }}:\nabla _x^S\delta {\bf{u}}} dv - \int\limits_{\varphi (\Omega )} {\delta {\bf{u}} \cdot {\bf{\bar b}}} dv - \int\limits_{\varphi ({\Gamma _t})} {\delta {\bf{u}} \cdot {\bf{\bar t}}} d\Gamma  = 0
\label{Eq6}
\end{equation}
By applying the Galerkin discretization, the approximate displacement field $\bf{u}$ can be expressed in the relation $\bf{u} = {\bf{N\bar u}}$, where $\bf{N}$ is a set of shape functions associated with nodal displacement vector ${\bf{\bar u}}$. By introducing the matrix $\bf{B}$, the gradient of the test function is defined as

\begin{equation}
\nabla _x^S\delta {\bf{u}} = {\bf{B}}\delta {\bf{\bar u}}
\label{Eq7}
\end{equation}
where $\bf{B}$ is the matrix of the derivatives of shape functions with respect to the current configuration and its definition for node $\textit{I}$ is described by

\begin{equation}
{{\bf{B}}^I} = \left[ {\begin{array}{*{20}{c}}
\vspace{2mm}	{\dfrac{{\partial {N_I}}}{{\partial x}}}&0\\
\vspace{2mm}	0&{\dfrac{{\partial {N_I}}}{{\partial y}}}\\
	{\dfrac{{\partial {N_I}}}{{\partial y}}}&{\dfrac{{\partial {N_I}}}{{\partial x}}}
	\end{array}} \right]
\label{Eq8}
\end{equation}\\
Substituting Eq. \ref{Eq7} into Eq. \ref{Eq8}, the finite element formulation can be expressed as

\begin{equation}
{\bf{R}}\left( {\bf{u}} \right) = \int\limits_{\varphi (\Omega )} {{{\bf{B}}^T}{\boldsymbol{\sigma }}} dv - \int\limits_{\varphi (\Omega )} {{{\bf{N}}^T}{\bf{\bar b}}} dv - \int\limits_{\varphi ({\Gamma _t})} {{{\bf{N}}^T}{\bf{\bar t}}} d\Gamma  = 0
\label{Eq9}
\end{equation}\\
where ${\bf{R}}\left( {\bf{u}} \right)$ is the residual vector, and terms in this equation involve the internal and external force, which are written by

\begin{equation}
{{\bf{F}}_{{\mathop{\rm int}} }} = \int\limits_{\varphi (\Omega )} {{{\bf{B}}^T}{\boldsymbol{\sigma }}} dv
\label{Eq10}
\end{equation}

\begin{equation}
{{\bf{F}}_{ext}} = \int\limits_{\varphi (\Omega )} {{{\bf{N}}^T}{\bf{\bar b}}} dv + \int\limits_{\varphi ({\Gamma _t})} {{{\bf{N}}^T}{\bf{\bar t}}} d\Gamma
\label{Eq11}
\end{equation}\\
Due to the presence of nonlinear effects on the weak form, the procedure of linearization needs to be applied into the Eq. \ref{Eq9}. This leads to the following definition, which is decomposed into the geometric and material part. In particular, the first term is related to the variation of the stress $\boldsymbol{\sigma }$, which depends on the material behavior causes material effects. The second one involves the current state of the stress, which is responsible for geometric effects.
%${C^E}$
\begin{equation}
\Delta {\bf{R}}({\bf{u}}) = {{\bf{K}}_{\tan }}\Delta {\bf{u}} = \int\limits_{\varphi (\Omega )} {\Delta {{\bf{B}}^T}{\boldsymbol{\sigma }}} dv + \int\limits_{\varphi (\Omega )} {{{\bf{B}}^T}\Delta {\boldsymbol{\sigma }}dv}
\label{Eq12}
\end{equation}
With the examination of hyper-elastic material, the constitutive relation is expressed as ${\boldsymbol{\sigma }} = {\bf{C}}^e:{\nabla ^S}{\bf{u}}$, where ${\bf{C}}^e$ is the constitutive elastic material tensor. Thus, the variational form of constitutive relation can be presented as $\Delta {\boldsymbol{\sigma }} = {\bf{C}}^e:{\nabla ^S}\Delta {\bf{u}}$, and the introduction of the matrix $\bf{G}$ which is the result of taking variation $\Delta {\bf{B}}$ in the relation $\Delta {\bf{B}} = {\bf{G}}\Delta {\bf{u}}$. Substituting these relations into Eq. \ref{Eq12}, it can be rewritten as

\begin{equation}
\Delta {\bf{R}}({\bf{u}}) = {{\bf{K}}_{\tan }}\Delta {\bf{u}} = \left( {\int\limits_{\varphi (\Omega )} {{{\bf{G}}^T}{\bf{M}}} {\bf{G}}dv + \int\limits_{\varphi (\Omega )} {{{\bf{B}}^T}{\bf{C}}^e{\bf{B}}dv} } \right)\Delta {\bf{u}}
\label{Eq13}
\end{equation}

The tangential stiffness matrix is separately performed by two parts as

\begin{equation}
{{\bf{K}}_{geo}} = \int\limits_{\varphi (\Omega )} {{{\bf{G}}^T}{\bf{MG}}dv{\rm{ }}}
\label{Eq14}
\end{equation}

\begin{equation}
{{\bf{K}}_{mat}} = \int\limits_{\varphi (\Omega )} {{{\bf{B}}^T}{\bf{C}}^e{\bf{B}}dv{\rm{ }}} {\rm{  }}
\label{Eq15}
\end{equation}
where $\bf{M}$ is the matrix of stress components defined in the form as

\begin{equation}
{\bf{M}} = \left[ {\begin{array}{*{20}{c}}
	\vspace{2mm}{{\sigma _{xx}}{{\bf{I}}_{2 \times 2}}}&{{\sigma _{xy}}{{\bf{I}}_{2 \times 2}}}\\
	{sym}&{{\sigma _{yy}}{{\bf{I}}_{2 \times 2}}}
	\end{array}} \right]
\label{Eq16}
\end{equation}

The matrix $\bf{G}$ containing the derivatives of the shape functions reported in \citep{Wriggers_2008} defined at node $\textit{I}$ is given by

\begin{equation}
{{\bf{G}}^I} = \left[ {\begin{array}{*{20}{c}}
	\vspace{2mm}{\dfrac{{\partial {N_I}}}{{\partial x}}}&0\\
	\vspace{2mm}0&{\dfrac{{\partial {N_I}}}{{\partial x}}}\\
	\vspace{2mm}{\dfrac{{\partial {N_I}}}{{\partial y}}}&0\\
	0&{\dfrac{{\partial {N_I}}}{{\partial y}}}
	\end{array}} \right]{\rm{  }}
\label{Eq17}
\end{equation}

In the scope of hyper-elastic materials undergoing the finite deformation analysis, the constitutive material and stress tensors are determined by taking the derivatives of the strain energy function $W$ with respect to the right Cauchy-Green tensor ${\bf{C}}$. In particular, the second Piola-Kirchhoff stress tensor ${\bf{S}}$ in terms of the initial configuration is expressed as
\begin{equation}
     {\bf{S}} = 2\dfrac{{\partial W}}{{\partial {\bf{C}}}}
     \label{Eq_stress_Piola}
\end{equation}
The corresponding constitutive elastic tensor is defined as
\begin{equation}
     {\bf{C}}^E = 4\dfrac{{\partial W }}{{\partial {\bf{C}}\partial {\bf{C}}}}
     \label{Eq_consti_mat}
\end{equation}
The transformation of stress and constitutive material tensor to the current configuration is given in the following relations\\
\begin{equation}
          {\boldsymbol{\sigma }} = \dfrac{1}{J}{\bf{FS}}{{\bf{F}}^T}
          \label{Eq_trans_stress}
\end{equation}\\
\begin{equation}
          C_{ijkl}^e = {F_{im}}{F_{jn}}{F_{kp}}{F_{lq}}C_{mnpq}^E
          \label{Eq_trans_C}
\end{equation}
where $J$ is the determinant of the deformation gradient $\bf{F}$. The definition of $\bf{F}$ in the current configuration is given by
\begin{equation}
          {\bf{F}} = \dfrac{{\partial {\bf{x}}}}{{\partial {\bf{X}}}} = {\left( {\dfrac{{\partial {\bf{X}}}}{{\partial {\bf{x}}}}} \right)^{ - 1}} = {\left( {\dfrac{{\partial ({\bf{x}} - {\bf{u}})}}{{\partial {\mathbf{x}}}}} \right)^{ - 1}} = {\left( {{\bf{I}} - \dfrac{{\partial {\bf{u}}}}{{\partial {\bf{x}}}}} \right)^{ - 1}}
          \label{Eq_deform_grad}
\end{equation}

For the Neo-Hookean compressible material, the strain energy is given in the following logarithmic function as reported in \citep{Belytschko_Liu_Moran_2013,Wriggers_2008,Rashetnia_2015}
\begin{equation}
          W \left( {\bf{C}} \right) = \frac{1}{2}\lambda {\left( {\log J} \right)^2} - \mu \log J + \dfrac{\mu}{2} \left( {trace\left( {\bf{C}} \right) - 3} \right)
          \label{Eq_strain_ener_compress}
\end{equation}
where $\lambda$ and $\mu$ are the Lamé parameters obtained from the relation of Young’s modulus and Poisson’s ratio.
The stress and the constitutive elastic tensor defined in the current configuration can be obtained as
\begin{equation}
          {\boldsymbol{\sigma }} = \frac{1}{J}\left( {\lambda \log J{\bf{I}} + \mu \left( {{\bf{F}}{{\bf{F}}^T} - {\bf{I}}} \right)} \right)
          \label{Eq_stress_compress}
\end{equation}
\begin{equation}
          C_{ijkl}^e = \lambda {\delta _{ij}}{\delta _{kl}} + \left( {\mu  - \lambda \log J} \right)\left( {{\delta _{ik}}{\delta _{jl}} + {\delta _{il}}{\delta _{kj}}} \right)
          \label{Eq_consti_compress}
\end{equation}

For incompressible materials, the value of the deformation gradient's determinant is almost unchanged and is equal to 1. Hence, the corresponding strain energy derived from Eq. \ref{Eq_strain_ener_compress} reads
\begin{equation}
          W\left( {\bf{C}} \right) = \dfrac{\mu}{2} \left( {trace\left( {\bf{C}} \right) - 3} \right)
          \label{Eq_strain_ener_incom}
\end{equation}
The corresponding second Piola-Kirchhoff stress is expressed by
\begin{equation}
          {\bf{S}} = 2\dfrac{{\partial W}}{{\partial {\bf{C}}}} - p{{\bf{C}}^{ - 1}} = \mu {\bf{I}} - p{{\bf{C}}^{ - 1}}
          \label{Eq_second_Piola_incompre}
\end{equation}
Being different from the description in Eq. \ref{Eq_stress_Piola}, the term of the hydrostatic pressure is added to the stress tensor owing to the effects of the incompressibility of materials. By assuming that the material is said to be plane stress, the Cauchy-Green tensor $\mathbf{C}$ could be expressed in the following form such that the incompressibility condition must be satisfied:
\begin{equation}
          {\mathbf{C}} = \left[ {\begin{array}{*{20}{c}}
      	{{C_{11}}}&{{C_{12}}}&0\\
	{{C_{21}}}&{{C_{22}}}&0\\
	0&0&{{C_{33}}}
	\end{array}} \right],{\rm{ with\,}}{{\rm{C}}_{33}} = \frac{{{t^2}}}{{{T^2}}}
          \label{Eq_constitutive_incompre_str}
\end{equation}
in which $t$ and $T$ are the thickness distributions measured in the deformed and undeformed configuration, respectively. The thickness arisen from the deformation is computed as
\begin{equation}
          {t^2} = \frac{{{T^2}}}{{\det \left( {{\bf{\bar C}}} \right)}}
          \label{Eq_thickness}
\end{equation}
The second Piola-Kirchhoff tensor given in Eq. \ref{Eq_second_Piola_incompre} can be rewritten as
\begin{equation}
         {\bf{S}} = \mu \left( {{\bf{\bar I}} - \det {{\left( {{\bf{\bar C}}} \right)}^{ - 1}}{{{\bf{\bar C}}}^{ - 1}}} \right)
         \label{Eq_Piola_incompre_re}
\end{equation}
where $\bf{\bar I}$ is the unit tensor $2\times 2$, and $\bf{\bar C}$ is the in-plane dilatation tensor. The corresponding Cauchy stress is given by
\begin{equation}
         {\boldsymbol{\sigma }} = \mu \frac{1}{J}{\bf{F}}\left( {{\bf{\bar I}} - \det {{\left( {{\bf{\bar C}}} \right)}^{ - 1}}{{{\bf{\bar C}}}^{ - 1}}} \right){{\bf{F}}^T}
         \label{Eq_Cauchy_incompre}
\end{equation}
The fourth-order tensor of the constitutive material in the incompressible conditions defined in the deformed configuration is
\begin{equation}
         C_{ijkl}^e = \mu \det {\left( {{\bf{\bar C}}} \right)^{ - 1}}\left( {2{\delta _{ij}}{\delta _{kl}} + \left( {{\delta _{ik}}{\delta _{jl}} + {\delta _{il}}{\delta _{kj}}} \right)} \right)
         \label{Eq_constitutive_incompre}
\end{equation}
%As regards the definition of material, a nearly incompressible Neo-Hookean material \citep{Belytschko_Liu_Moran_2013,Wriggers_2008} through Cauchy stress tensor is employed in the finite elastic deformation. The formulation of the stresses is given by
%
%\begin{equation}
%{\boldsymbol{\sigma }} = \left[ {\begin{array}{*{20}{c}}
%\vspace{2mm}	{{\sigma _{xx}}}&{{\sigma _{xy}}}\\
%	{{\sigma _{yx}}}&{{\sigma _{yy}}}
%	\end{array}} \right] = \dfrac{1}{J}\left( {\lambda \ln J{\bf{I}} + \mu \left( {{\bf{F}}{{\bf{F}}^T} - {\bf{I}}} \right)} \right)
%\label{Eq18}
%\end{equation}
%where $\lambda$ and $\mu$ are Lamé constants of the linearized theory which are converted from Young’s modulus and Poisson’s ratio, and $J$ is the determinant of the deformation gradient $\mathbf{F}$ which is defined on the current configuration is given by
%\begin{equation}
%{\bf{F}} = \frac{{\partial {\bf{x}}}}{{\partial {\bf{X}}}} = {\left( {\frac{{\partial {\bf{X}}}}{{\partial {\bf{x}}}}} \right)^{ - 1}} = {\left( {\frac{{\partial ({\bf{x}} - {\bf{u}})}}{{\partial {\mathbf{x}}}}} \right)^{ - 1}} = {\left( {{\bf{I}} - \frac{{\partial {\bf{u}}}}{{\partial {\bf{x}}}}} \right)^{ - 1}}
%\label{Eq3}
%\end{equation}

\subsection{Large deformation formulation in XFEM}
\label{sub_sec 2.2}

XFEM is a well-known technique in modeling discontinuities without remeshing. Being proposed by Belytschko \textit{et al} \citep{Belytschko_Black_1999, Moes_Dolbow_1999}, the method mathematically described the effects of discontinuities by additional terms, namely enrichment functions. In modeling crack problems, two types of enrichments are utilized, including the Heaviside function to express the jump in displacements and the asymptotic functions to capture the singularity of stresses near the crack tip. The approximation displacement fields in the crack problem is given by

\begin{equation}
{\bf{u}}({\bf{x}}) = \sum\limits_{i \in I} {{N_i}({\bf{x}}){{{\bf{\bar u}}}_i}}  + \sum\limits_{j \in J} {{N_j}({\bf{x}})H({\bf{x}}){{{\bf{\bar a}}}_j}}  + \sum\limits_{k \in K} {{N_k}({\bf{x}})\sum\limits_l {{A_l}({\bf{x}})} {\bf{\bar b}}_k^l}
\label{Eq19}
\end{equation}
where ${N_i}({\bf{x}})$ and ${{\bf{\bar u}}_i}$ are the shape functions and the nodal displacement of node $\textit{i}$, and ${{\bf{\bar a}}_i}$ is the nodal enriched degree of freedom associated with the Heaviside function, while ${\bf{\bar b}}_k^l$ is the nodal enriched degree of freedom corresponding to node $\textit{k}$ associated with the asymptotic enrichment. The Heaviside function is defined by a step function

\begin{equation}
H({\bf{x}}) = \left\{ {\begin{array}{*{20}{c}}
	1&{{\rm{if (}}{\bf{x}} - {{\bf{x}}^*}) \cdot {\bf{n}} > 0}\\
	{ - 1}&{{\rm{other}}}
	\end{array}} \right.
\label{Eq20}
\end{equation}
where ${\bf{x}}^*$ is the closet point on the crack to point $\bf{x}$, and $\bf{n}$ is the normal vector to crack at ${\bf{x}}^*$. \\

The enrichment functions related to the asymptotic field near the crack tip are usually examined under the use of four branch functions \citep{Belytschko_Black_1999,Moes_Dolbow_1999} taken from the linear elastic fracture mechanics. Several studies on finite strains around the crack tip were analyzed in the case of plane stress shown in \citep{Philippe_Wolfgang_1994_1} and plane strain \citep{Knowles_Sternberg_1973}, based on the use of Neo-Hookean strain energy. Legrain et al \citep{Legrain_Moes_Verron_2005} then proposed the enrichment function to deal with this problem appropriate to the framework of XFEM. It is given by

\begin{equation}
{A}({\bf{x}}) = \left\{ {{r^{1/2}}\sin (\theta /2)} \right\}
\label{Eq21}
\end{equation}
where $r$ and $\theta$ are the polar coordinates with respect to the crack tip.\\

\begin{figure}
	\centering
	\includegraphics[scale=1]{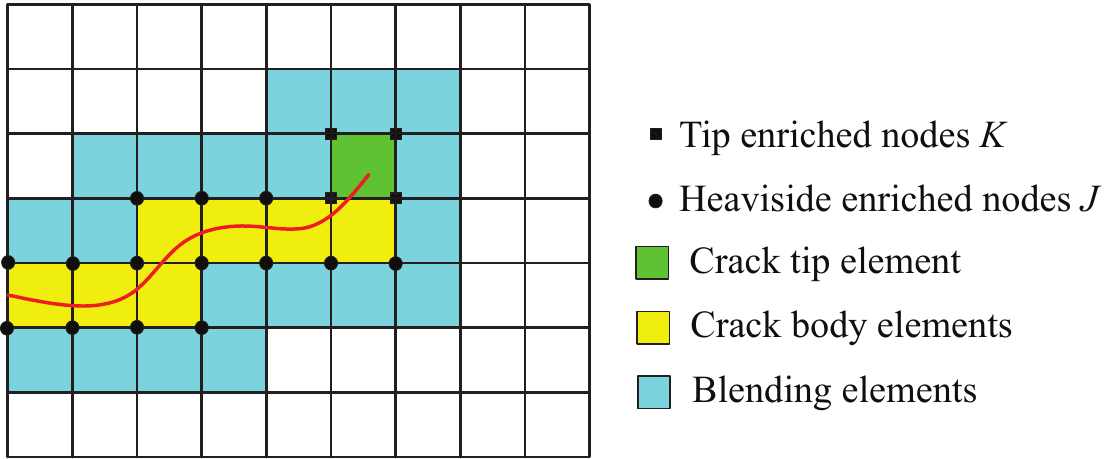}
	\caption{The crack path in a domain with fully enriched elements and blending elements.}
	\label{F2}
\end{figure}

An important feature in the XFEM analysis is that a certain number of elements are partially enriched, whose are called ``blending element" as depicted in Fig. \ref{F2}. Therefore, the sum of basis functions of these elements does not satisfy the property of the Partition of Unity. As a result, the enrichment functions cannot be exactly reproduced the approximation fields. This also causes poor convergence rates. A solution proposed by Fries \citep{Fries_2008} is used to treat the negative effects of blending elements, in particular, a modified approach relied on Ramp function $R\left( {\bf{x}} \right)$ which is established from the shape functions of possible enriched nodes $K$ belonging to counterparts of fully enriched elements for the crack tip:

\begin{equation}
R\left( {\bf{x}} \right) = \sum\limits_{k \in K} {{N_k}\left( {\bf{x}} \right)} {\rm{ }}
\label{Eq22}
\end{equation}

The modification of the asymptotic enrichment functions are defined as

\begin{equation}
{{\bf{A}}^{\bmod }}\left( {\bf{x}} \right) = {\bf{A}}\left( {\bf{x}} \right)R\left( {\bf{x}} \right)
\label{Eq23}
\end{equation}
The definition of the modified enrichment functions is zero in the standard elements and continuously varies in between elements.\\

For the implementation of XFEM, enrichment functions are often expressed in the shifted relation, in which take the value of zero at nodal points. By associating this with the modified enrichment functions, the approximation displacements are rewritten as

\begin{equation}
{\bf{u}}({\bf{x}}) = \sum\limits_{i \in I} {{N_i}({\bf{x}}){{{\bf{\bar u}}}_i}}  + \sum\limits_{j \in J} {{N_j}({\bf{x}})\bar H({\bf{x}}){{{\bf{\bar a}}}_j}}  + \sum\limits_{k \in K} {{N_k}({\bf{x}})\sum\limits_l {{{\bar A}_l}({\bf{x}})R({\bf{x}})} {\bf{\bar b}}_k^l}
\label{Eq24}
\end{equation}
where $\bar H({\bf{x}}) = H({\bf{x}}) - H({{\bf{x}}_j})$, and ${\bar A_l}({\bf{x}}) = {A_l}({\bf{x}}) - {A_l}({{\bf{x}}_k})$.\\

\begin{figure}
	\centering
	\includegraphics[scale=1]{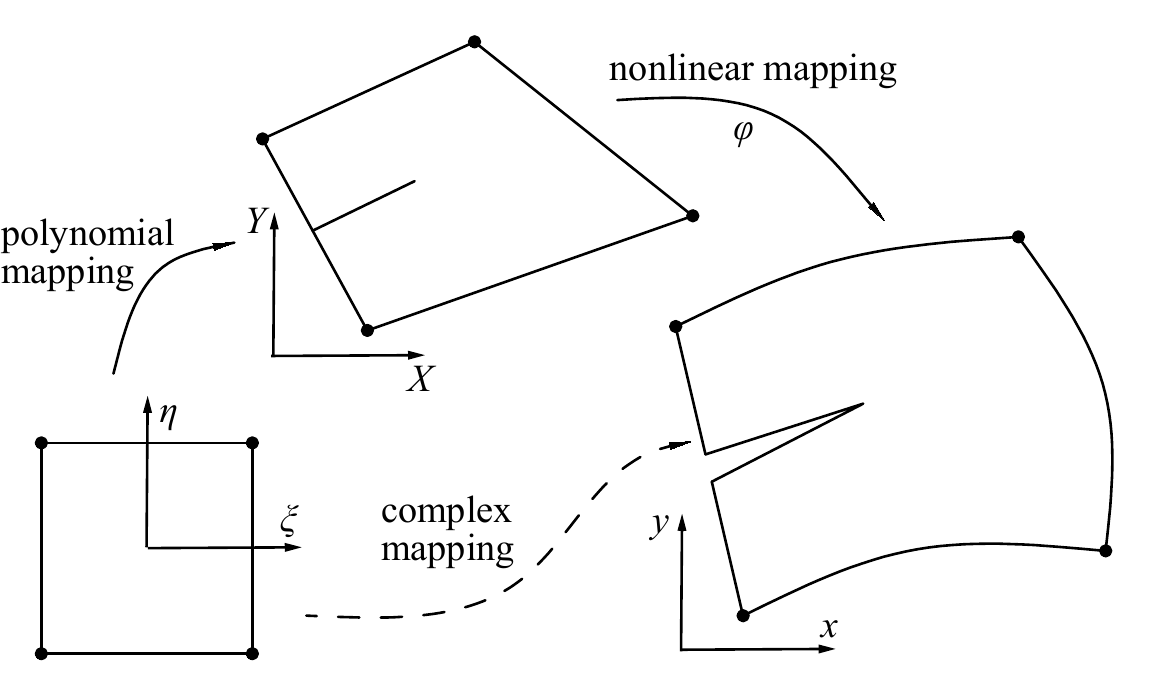}
	\caption{An illustration of mappings of an enriched element.}
	\label{F3}
\end{figure}

As mentioned earlier, the nonlinear crack problem by using XFEM is implemented in the updated Lagrangian description. The mapping between physical (ordinary) elements and their parent counterparts is polynomial, and is straightforward for any configurations. However, the mapping from the parent elements in the natural coordinates to their enriched elements in the current configuration is very complicated or even is impossible because the enrichment functions adapted to this element are with respect to the polar coordinates including $r$ and $\theta$. In general, the value of the angle $\theta$ is not guaranteed under the process of mapping. Thus, enrichment functions must be defined on the initial configuration, and then a transformation $\varphi$ is employed in the mapping from the initial (undeformed) configuration to the current (deformed) configuration, reported in \citep{Broumand_Khoei_2013} as illustrated in Fig. \ref{F3}. It is noted that the second mapping stage does not remain polynomial due to the effect of the enrichment function.\\

\begin{figure}
	\centering
	\includegraphics[scale=1]{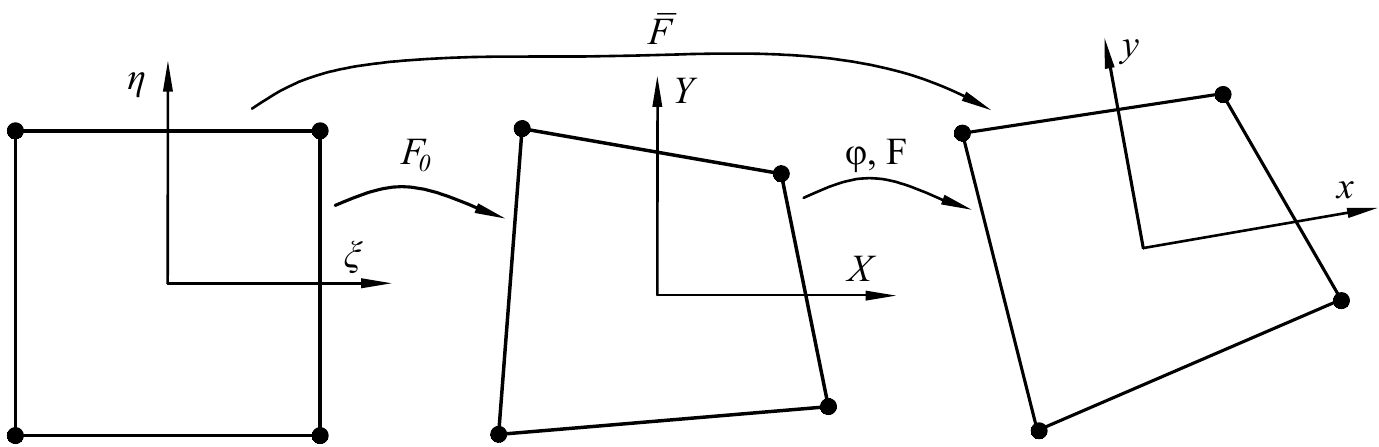}
	\caption{The process of mapping of an ordinary element in the large deformation XFEM.}
	\label{F4}
\end{figure}

For the mapping of ordinary elements, the process of mapping between two different configurations is performed in the flowchart as shown in Fig. \ref{F4}. In particular, the mapping from the reference element to the undeformed configuration is denoted by ${F_0}$, and the other manner from undeformed to deformed configuration is denoted by $F$, also called the deformation gradient. The reference element can be mapped directly to the deformed configuration, denoted by $\bar F$.\\

The definitions of the mapping of ${F_0}$ and $F$ are indirectly determined through the following computation of the Jacobian matrices

\begin{equation}
{{\bf{J}}_0} = \left[ {\begin{array}{*{20}{c}}
	\vspace{3mm}{\dfrac{{\partial X}}{{\partial \xi }}}&{\dfrac{{\partial Y}}{{\partial \xi }}}\\
	{\dfrac{{\partial X}}{{\partial \eta }}}&{\dfrac{{\partial Y}}{{\partial \eta }}}
	\end{array}} \right]{\rm{ }}
\label{Eq25}
\end{equation}
\begin{equation}
{{\bf{J}}} = \left[ {\begin{array}{*{20}{c}}
	\vspace{3mm}{\dfrac{{\partial x}}{{\partial \xi }}}&{\dfrac{{\partial y}}{{\partial \xi }}}\\
	{\dfrac{{\partial x}}{{\partial \eta }}}&{\dfrac{{\partial y}}{{\partial \eta }}}
	\end{array}} \right]{\rm{ }}
\label{Eq26}
\end{equation}
\begin{equation}
{{\bf{F}}_0}{\rm{  =  }}{\bf{J}}_0^{ - 1}
\label{Eq27}
\end{equation}
\begin{equation}
{\bf{\bar F}}{\rm{  =  }}{{\bf{J}}^{ - 1}}
\label{Eq28}
\end{equation}

The term of $\bf{F}$ is the deformation gradient, which is given in the Eq. \ref{F3} is also defined as
\begin{equation}
{\bf{F}} = \left[ {\begin{array}{*{20}{c}}
	\vspace{3mm}{\dfrac{{\partial x}}{{\partial X}}}&{\dfrac{{\partial x}}{{\partial Y}}}\\
	{\dfrac{{\partial y}}{{\partial X}}}&{\dfrac{{\partial y}}{{\partial Y}}}
	\end{array}} \right]{\rm{ = }}\left[ {\begin{array}{*{20}{c}}
	\vspace{3mm}{{J^{11}}F_0^{11} + {J^{21}}F_0^{12}}&{{J^{11}}F_0^{21} + {J^{21}}F_0^{22}}\\
	{{J^{12}}F_0^{11} + {J^{22}}F_0^{12}}&{{J^{12}}F_0^{21} + {J^{22}}F_0^{22}}
	\end{array}} \right]
\label{Eq29}
\end{equation}

In the same way to define the approximation displacement as shown in Eq. \ref{Eq24}, the current position in the framework of XFEM can be given by

\begin{equation}
{\bf{x}} = \sum\limits_{i \in I} {{N_i}(\xi ,\eta ){{{\bf{\bar x}}}_i}}  + \sum\limits_{j \in J} {{N_j}(\xi ,\eta )\bar H({\bf{x}}){{{\bf{\bar a}}}_j}}  + \sum\limits_{k \in K} {{N_k}(\xi ,\eta )\sum\limits_l {{{\bar A}_l}({\bf{X}})R({\bf{X}})} {\bf{\bar b}}_k^l}
\label{Eq30}
\end{equation}
Based on the definitions of Jacobian matrix as reported in Eq. \ref{Eq25} and Eq. \ref{Eq26}, their corresponding calculations in XFEM are presented in \ref{appendixA}. The construction of the matrices $\bf{B}$, and $\bf{G}$ shown in Eq. \ref{Eq20} and Eq. \ref{Eq21} can be derived from the context of XFEM in relation with the terms of enrichments. Therefore, the components of these matrices are comprised of three different parts, including the standard FEM, Heaviside enrichment and asymptotic enrichments. More detail for their formulations can be found in \ref{appendixB}.\\

For the numerical treatment of the nonlinear equations, the Newton-Raphson’s method is widely employed to obtain solutions according to the iterative procedure to find the displacement increment $\Delta {\bf{u}}$. As mentioned in Eq. \ref{Eq13}, the updated state of the residual force corresponding to ${k^{th}}$ iteration of the Newton-Raphson method can be defined as

\begin{equation}
{{\bf{R}}_{k + 1}} \simeq {{\bf{R}}_k} + \Delta {\bf{R}}
\label{Eq31}
\end{equation}
The vector $\bf{R}_{k + 1}$ can be described in the following expression ${{\bf{R}}_{k + 1}} = {\bf{R}}({{\bf{u}}_k} + \Delta {\bf{u}})$, which is equal to ``$\bf{0}$" as defined in Eq. \ref{Eq9}. Therefore, the formulation to get the incremental displacements in Eq. \ref{Eq13} can be rewritten at an iterative state as

\begin{equation}
{{\bf{K}}_{\tan }}\Delta {\bf{u}} =  - {{\bf{R}}_k}
\label{Eq32}
\end{equation}
For the models with many load steps, the convergent solution is set up as a starting value of the unknown for the next load stage. An algorithm for the framework of the nonlinear XFEM implementation is reported in Table \ref{tab:Algorithm}.

\begin{table}[!htb]
	\centering
	\begin{tabular}{|l|p{5.9cm}|}
		\hline
		\textbf{Algorithm} Algorithm for the Newton-Raphson iteration\\
		\hline
		\quad 1. Make an initial displacement $\mathbf{u}_0$, normally $\mathbf{u}_0 = \mathbf{0}$\\
		\quad 2. \textbf{Loop} over each incremental load\\
		\qquad a. \textbf{Loop} until the converged solution is achieved, in particular $\left\| {{{\bf{R}}_k}} \right\| > tol$, \\
		\quad\qquad where $\mathbf{R}_k$ is obtained in Eq. \ref{Eq9}, and $k$ is iteration steps.\\
		\quad\qquad i.  Calculate the tangential stiffness matrix ${{\bf{K}}_{tan }}$, internal force ${{\bf{K}}_{int}}$, and \\
		\qquad\qquad external force ${{\bf{K}}_{ext}}$ at the current displacement $\mathbf{u}_k$, mentioned in Eq. \ref{Eq13},\\
		\qquad\qquad Eq. \ref{Eq10}, and Eq. \ref{Eq11}, respectively. \\
		\quad\qquad ii. Calculate the residual force ${{\bf{R}}_k} = {{\bf{F}}_{int}} - {{\bf{F}}_{ext}}$.\\
		\quad\qquad iii.Calculate the incremental displacement $\Delta {\bf{u}} = {({{\bf{K}}_{\tan }})^{ - 1}}( - {{\bf{R}}_k})$.\\
		\quad\qquad iv. Update the displacement ${{\bf{u}}_{k + 1}} = {{\bf{u}}_k} + \Delta {\bf{u}}$.\\
		\qquad b. \textbf{End Loop}.\\
		\qquad c. \textbf{Increase load}.\\
		\quad 3. \textbf{End Loop}\\ \hline
	\end{tabular}
	\caption{Algorithm for solving convergent solution in the nonlinear XFEM}
	\label{tab:Algorithm}
\end{table}	

\section{Large deformation formulation based on XFEM over polygonal meshes}
\label{sec 3}
\subsection{Polygonal finite elements}
\label{sub_sec 3.1}

The concept of polygonal elements is extended from standard triangular and quadrilateral elements. The number of element vertices is arbitrary, and the shape can even be concave. The early contribution in this field is the development of convex multi-sided elements using rational basis functions proposed by Wachspress \citep{Wachspress_1979}. Laplace shape functions were accordingly constructed by the use of natural neighbor interpolation for arbitrary convex polygons introduced by Sibson\citep{Sibson_2012}. Floater \citep{Floater_2014,Floater_2016} relied on the application of barycentric coordinates to develop Wachspress and mean value basis functions which are commonly employed in computer graphics and numerical modeling. Hormann \citep{Hormann_2006} successively exploded the mean value coordinates to be well-defined for arbitrary planar and even non-convex polygons. Besides that, the mean value coordinates are suitable for the larger deformation analysis when a number of possible elements are distorted. In the present work, the use of mean value shape functions is taken into the conforming approximations on polygonal finite elements. Following that, a modification regarding gradient correction of polygonal shape functions is build up to enhance their accurate computation.\\

\begin{figure}[!httb]
	\centering
	\includegraphics[scale=1]{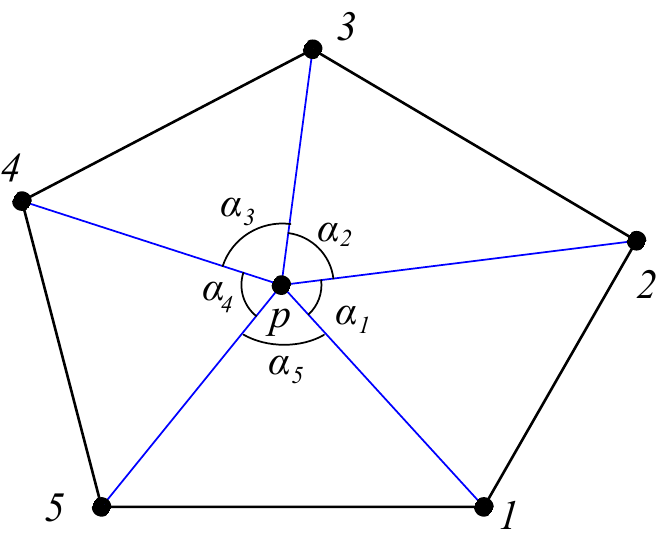}
	\caption{Notations for mean value coordinates of a pentagon at point $\mathit{p}$.}
	\label{Fig.5}
\end{figure}

Basis functions of the mean value concept are constructed from non-negative weights on polygons. The polygonal shape functions of a $n$-gons domain are expressed in general form as

\begin{equation}
	{N_i}\left( {\bf{x}} \right) = \frac{{{w_i}\left( {\bf{x}} \right)}}{{\sum\limits_j^n {{w_j}\left( {\bf{x}} \right)} }},\quad {\rm{ }}{w_i}\left( {\bf{x}} \right) = \frac{{\tan ({\alpha _{i - 1}}/2) + \tan ({\alpha _i}/2)}}{{\left\| {{{\bf{v}}_i} - {\bf{x}}} \right\|}}{\rm{ }}
	\end{equation}
	where ${\bf{v}}_i$ are nodal positions, and $\bf{x}$ is the position of point $p$ in the domain. The value of angles ${\alpha _i}({\bf{x}})$ is defined such that $0 < {\alpha _i} < \pi $. Fig. \ref{Fig.5} describes notations for mean value coordinates of a pentagon.

The mean value shape functions satisfy the interpolation requirements, including:\\
\begin{itemize}
	\item A partition of unity, and non-negative character: $\sum\limits_{i = 1}^n {{N_i}\left( {\bf{x}} \right)}  = 1$, $0 \le {N_i}\left( {\bf{x}} \right) \le 1$.
	\item Interpolate nodal data: ${N_i}\left( {{{\bf{x}}_j}} \right) = {\delta _{ij}}$, where ${\delta _{ij}}$ is the Kronecker-delta.
	\item Linear completeness: ${\bf{x}} = \sum\limits_{i = 1}^n {{N_i}\left( {\bf{x}} \right){{\bf{x}}_i}}$.
\end{itemize}

\begin{figure}[!httb]
	\centering
	\subfigure[]{
		\includegraphics[scale=1]{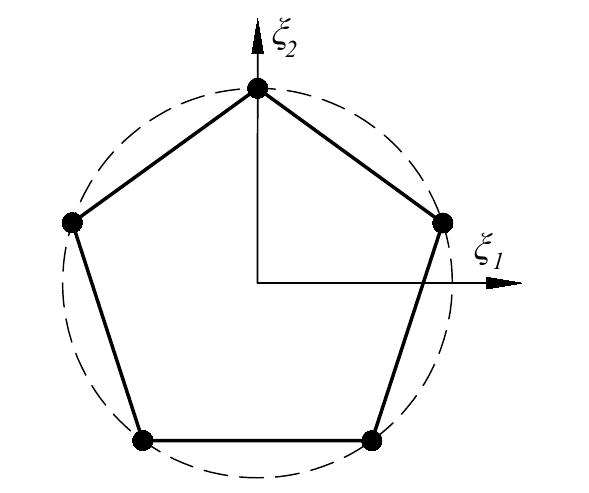}
		\label{Fig.6a}
	} \hspace*{1em}
	\subfigure[]{
		\includegraphics[scale=1]{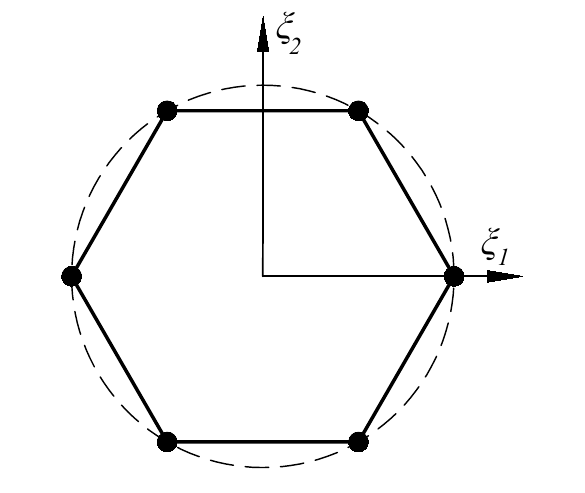}
		\label{Fig.6b}
	} \hspace*{0em}\\
	\caption{Reference elements: (a) pentagon, (b) hexagon.}
	\label{Fig.6}
\end{figure}

Conforming approximations on a polygonal element in physical space can be implemented from the definition of shape functions directly defined on physical coordinates. A common way in numerical analysis is the employment of an iso-parametric mapping from the corresponding reference element. For the convenience in implementation, all of vertices of canonical elements lie on the same circumcircle and are arranged in the anticlockwise order. An example of pentagon and hexagon reference elements is shown in Fig. \ref{Fig.6}.\\

\begin{figure}
	\centering
	\includegraphics[scale=1]{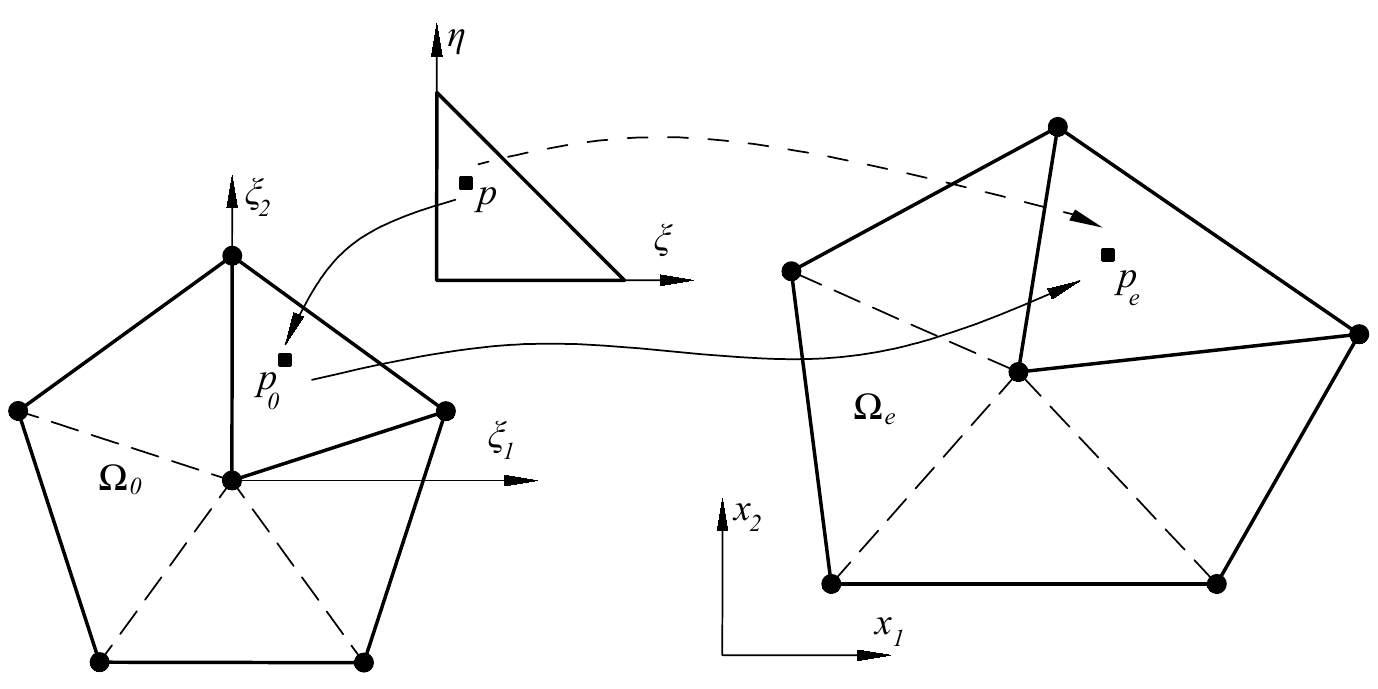}
	\caption{Numerical integration based on partition schemes.}
	\label{Fig.7}
\end{figure}

To evaluate the numerical integration over polygonal elements, a popular technique is the use of two levels of affine mapping, reported in \citep{Sukumar_Tabarraei_2004,Sukumar_Malsch_2006,Tabarraei_Sukumar_2008}, which can be straightforward to apply Gaussian quadrature. A depiction of the method is shown in Fig. \ref{Fig.7}. The first level is the mapping from a polygonal domain in physical coordinate to a corresponding reference domain. The reference domain then is divided into sub-triangles, and each sub-triangle is mapped into a reference triangle. Taking the integration of a typical function over a physically function $f$ over a polygonal domain is expressed as
\begin{equation}
    \int_{{\Omega _e}} {fd\Omega }  = \int_{{\Omega _0}} {f\left| {{{\bf{J}}_2}} \right|d\Omega }  = \sum\limits_{i = 1}^n {\int_0^1 {\int_0^{1 - \xi } {f\left| {{{\bf{J}}_1^i}} \right|\left| {{{\bf{J}}_2}} \right|d\xi d\eta } } }
    \label{Eq34}
\end{equation}
Alternatively, the numerical integration can be executed by a direct partition of the physical element into sub-triangles to apply quadrature rules. The integral process is simplified in one level as shown in Fig. \ref{Fig.7}, and Eq. \ref{Eq34} is thus written as
	\begin{equation}
	\int_{{\Omega _e}} {fd\Omega }  = \sum\limits_{i = 1}^n {\int_0^1 {\int_0^{1 - \xi } {f\left| {{\bf{J}}_1^i} \right|d\xi d\eta } } }
	\end{equation}

It is observed that the basis functions mentioned earlier are of much benefit to conforming to arbitrary polygonal finite elements when the interpolation conditions are completely satisfied. Nevertheless, the numerical integration of gradient fields on such domains causes errors and instabilities in convergent results because of their non-polynomial property. A commonly recent solution introduced by Talischi et al \citep{Talischi_2015,Chi_2016} is the use of gradient correction to generate polynomial consistency which definitely enhances the accuracy of polygonal finite element methods. This modification is briefly expressed as follows. For functions $N$ defined on a domain $\Omega_e$, the gradient correction of $N$ denoted by ${\nabla _{\Omega_e,k}}N$ that is close to the original form $\nabla N$ relatively requires the surface and line integration. The formulation of the correction gradient is given by
\begin{equation}
	{\nabla _{\Omega_e,k}}N = \nabla N + \frac{1}{{\left| \Omega_e \right|}}\left( {\oint_{\partial \Omega_e} {N{\bf{n}}ds}  - \int_{\Omega_e} {\nabla Nd{\bf{x}}} } \right)
\end{equation}
where the subscript “$k$” denotes an elemental order.\\

\begin{figure}
	\centering
	\subfigure[]{
		\includegraphics[scale=1]{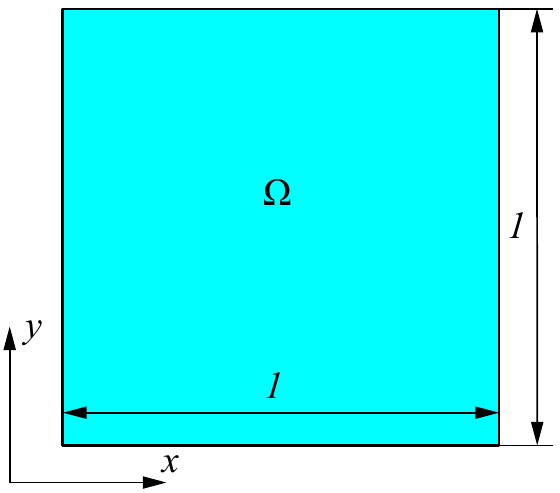}
		\label{Fig.patch}
	} \hspace*{3em}
	\subfigure[]{
		\includegraphics[scale=1]{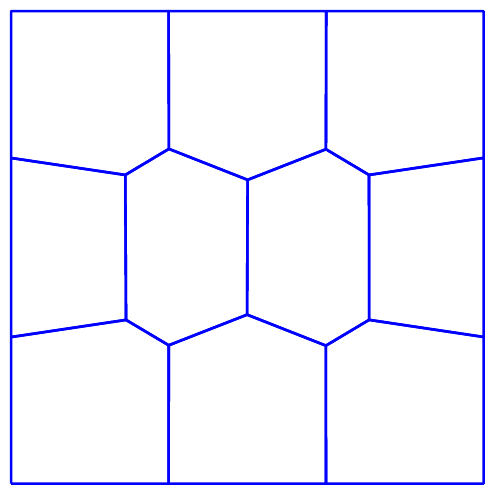}
		\label{Fig.patch_mesh}
	} \hspace*{0em}\\
	\caption{(a) Geometry of the patch test, (b) a typical polygonal mesh.}
	\label{Fig.patch_test}
\end{figure}

An investigation for validating the computational efficiency of gradient correction technique for basis functions is employed in a patch test problem bounded in the domain $\Omega  = {\left( {0,{\rm{ }}1} \right)^2}$ as shown in Fig. \ref{Fig.patch}. The analytical solution for displacements is given by ${u_x}(x,y) = 2x$, and ${u_y}(x,y) =-0.5y$, along with material parameters $\mu  = \kappa  = 1$. The verification is carried out on meshes of linear polygonal elements $(k=1)$ as typically presented in Fig. \ref{Fig.patch_mesh} to provide capacity of the correction scheme for mean value shape functions. Relative errors in the displacement field are studied through evaluation of ${L^2}$-error $\left\| {{\bf{u}} - {{\bf{u}}_h}} \right\|$ and ${H^1}$-energy norm $\left\| {\nabla {\bf{u}} - \nabla {{\bf{u}}_h}} \right\|$. Resulting data presented in Table. \ref{Tab2} proves that the correction gradient has significant effects on mean value basis functions. Obviously, the correction approach raises remarkable accuracy when the patch test is exactly passed.

\begin{table}
	\centering	
	\caption{Results of relative errors in ${L}^2$ and ${H}^1$ norms for the patch test.}
	\begin{tabular}{|c|c|c|c|c|}
		\hline
		Number of & \multicolumn{2}{c|}{${L}^2$-error} & \multicolumn{2}{c|}{${H}^1$-energy norm} \\
		 \cline{2-5}elements & Correction & Without correction & Correction & Without correction\\
		 \hline
		 10    & 3.45E-16 & 3.09E-04 & 1.71E-15 & 1.40E-03 \\
		 50    & 5.01E-16 & 1.11E-04 & 3.32E-15 & 8.86E-04 \\
		 250   & 7.12E-16 & 6.02E-05 & 7.05E-15 & 1.09E-03 \\
		 \hline
		
	\end{tabular}
	\label{Tab2}
\end{table}

\subsection{Local mesh refinements with hanging nodes}
\label{sub_sec 3.2}

\begin{figure}
	\centering
	\includegraphics[scale=1]{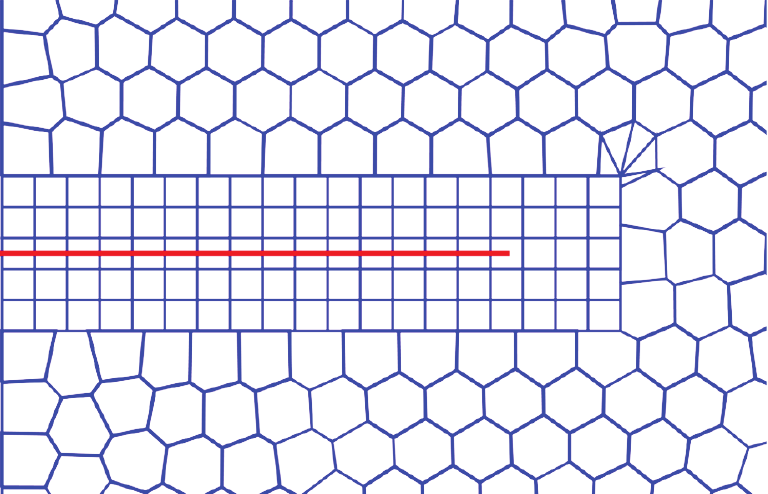}
	\caption{A mesh with hanging nodes.}
	\label{Fig.8}
\end{figure}

The applications of polygonal finite element into XFEM are certainly beneficial to enhance the computational efficiency thanks to the high order of shape functions of such elements, and the flexibility in the mesh generation. In the finite deformation analysis, the presence of displacement jumps where possibly causes mesh distortion is necessary for the treatment. Meshes generated by polygonal elements are unstructured, therefore regions with significantly large deformations are difficult to control the mesh distortion phenomenon. A proposed solution is the addition of structured meshes around the area where attending the discontinuities. The practice of using rectangle-shaped elements around the crack interface is useful to reduce mesh distortions when the effect of aspect-ratio distortions at such split elements is on mesh quality. Due to the differences from the size and the type of elements along the boundaries of this area, hanging nodes are produced, which lead to incompatibilities in the mesh. A typical example of the local refinement of a crack illustrates the presence of hanging nodes shown in Fig. \ref{Fig.8}.\\

\begin{figure}
	\centering
	\includegraphics[scale=1]{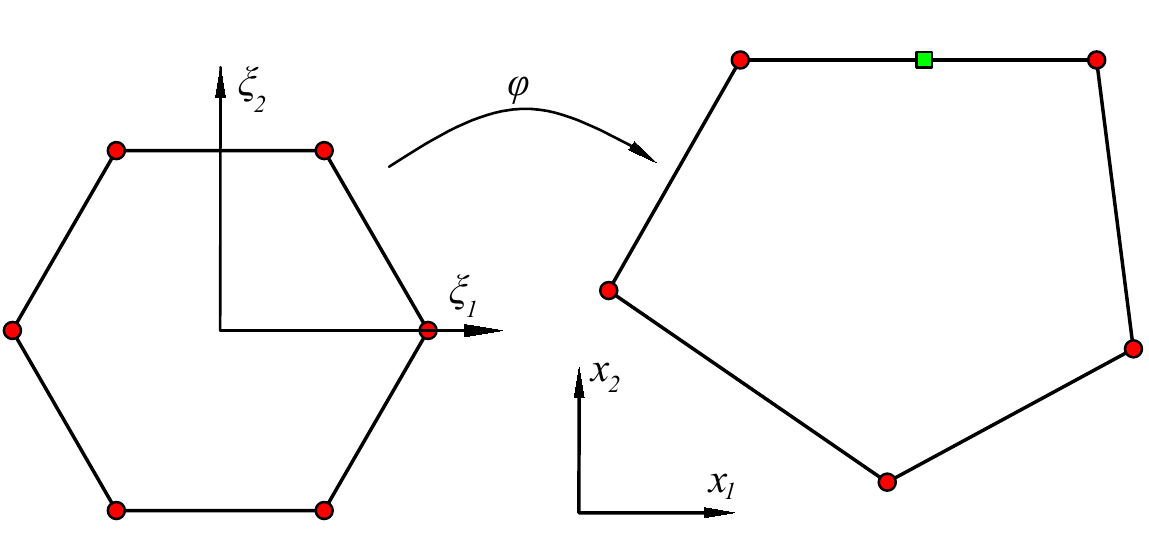}
	\caption{Mapping from a canonical hexagon to a pentagon with one hanging node.}
	\label{Fig.9}
\end{figure}

The treatment of hanging nodes can be made by two common approaches which are relied on the consideration of degrees of freedom (DOFs). In particular, the first approach is the case of hanging nodes having no DOFs and non-conforming shape functions. The numerical solutions are obtained by considering constrained approximation via the Lagrangian multipliers, which result in complicated algorithms. In the second one, the hanging nodes are associated with DOFs, and therefore the construction of shape functions is required to conform to finite element space. Gupta \citep{Gupta_1978} developed shape functions based on bilinear basis functions on quadrilateral elements. However, there is the limit to the number of hanging nodes on a certain element edge. The application of polygonal finite elements to conform to quadtree meshes has been interesting in adaptive computations \citep{Tabarraei_Sukumar_2008,Tabarraei_2007}. With the great flexibility of the shape functions, there is no distinction between hanging nodes and regular nodes. Therefore, the use of shape functions over polygonal elements are well-defined on elements with the presence of hanging nodes. Fig. \ref{Fig.9} shows the generation of a hanging node from a pentagon which is accordingly mapped into the reference hexagon.

\section{Numerical examples}
\label{sec 4}

In this section, several examples are investigated to demonstrate the computational efficiency of the present approach. The Neo-Hookean models mentioned in Eq. \ref{Eq_strain_ener_compress} and Eq. \ref{Eq_strain_ener_incom} for compressible and incompressible hyper-elastic materials are employed in modeling the nonlinear fracture problems. The first examples are made to verify the reliability of the present numerical technique, and prove nonlinearities of the structural mechanics in comparison with previous study counterparts regarding both numerical and analytical results. The advantages of applying polygonal finite elements into the XFEM are exhibited in the next examinations, whose benchmarks are included with a mixed crack or complicated geometries. The Lamé parameters defined in the material models in Section \ref{sub_sec 2.1} are converted from the relations of Young's modulus and Poisson's ratio found in \citep{Wriggers_2008} as
\begin{equation}
    E = \frac{{(3\lambda  + 2\mu )\mu }}{{\lambda  + \mu }}
    \label{Eq35}
\end{equation}

\begin{equation}
   \nu  = \frac{\lambda }{{2(\lambda  + \mu )}}
   \label{Eq36}
\end{equation}
Note that benchmarks with incompressible elastic materials for the numerical analysis here are assumed to be in the plane stress state. The convergence tolerance is set to be $tol = 6 \times {10^{ - 3}}$.

%A key assessment in the fracture mechanics is the determination of the stress intensity factors (SIFs). The common way to evaluate SIFs is the use of the $J$-integral, which is specified in the line integral to compute the energy release rate for a crack. The $J$-integral defined in the deformed configuration reported in \citep{Legrain_Moes_Verron_2005} takes the following form
%
%\begin{equation}
%    J = \int_\gamma  {\left( {w{n_1} - {\bf{n}} \cdot {\boldsymbol{\sigma}} \cdot \frac{{\partial {\bf{u}}}}{{\partial {x_1}}}} \right)d\gamma }
%    \label{Eq37}
%\end{equation}
%where $w$ is the strain energy, $\bf{n}$ is the unit normal outward vector, and $\gamma$ is the integration contour. A depiction of the $J$-integral is shown in Fig. \ref{Fig.10}.
%
%\begin{figure}
%	\centering
%	\includegraphics[scale=1]{F10_integral_domain.pdf}
%	\caption{The contour of the $J$-integral.}
%	\label{Fig.10}
%\end{figure}

\subsection{A single edge crack specimen in a square plate}
\label{sub_sec 4.1}

\begin{figure}
	\centering
	\includegraphics[scale=1]{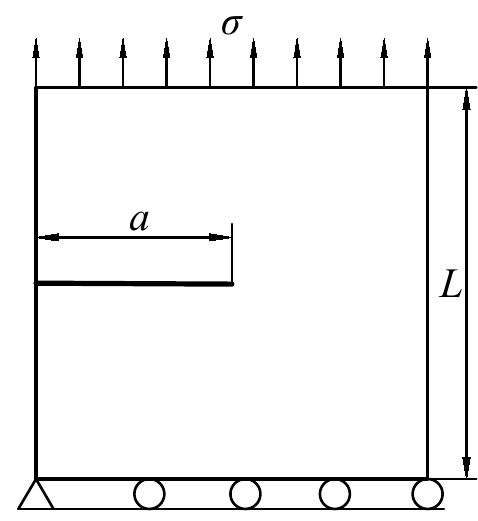}
	\caption{Geometry and boundary conditions of the single edge crack problem.}
	\label{Fig.11}
\end{figure}

\begin{figure}[!httb]
	\centering
	\subfigure[]{
		\includegraphics[scale=1]{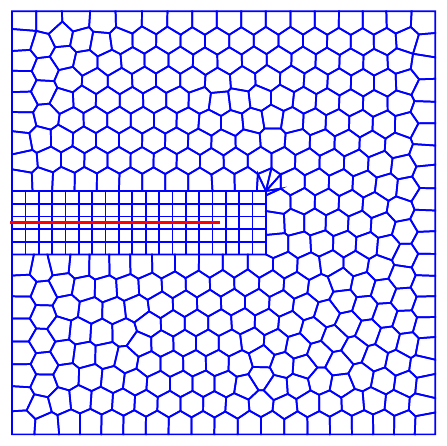}
		\label{Fig.11a}
	} \hspace*{0em}
	\subfigure[]{
		\includegraphics[scale=1]{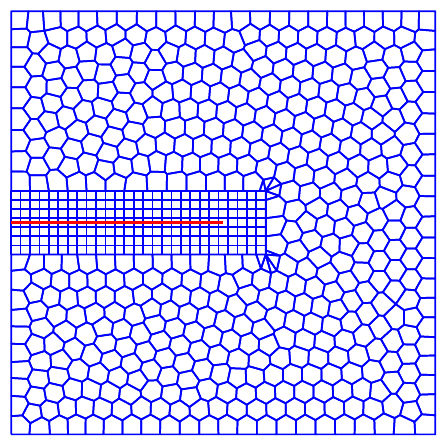}
		\label{Fig.11b}
	} \hspace*{0em}
    \subfigure[]{
    	\includegraphics[scale=1]{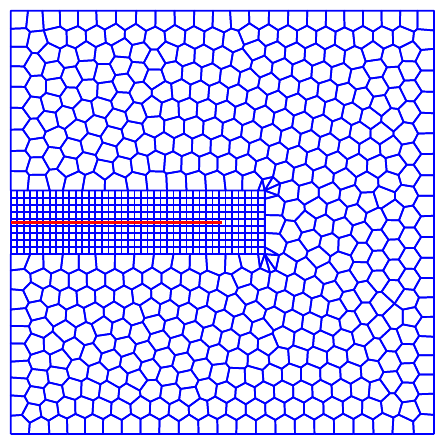}
    	\label{Fig.11c}
    } \hspace*{0em}
	
	\caption{Polygonal mesh discretization of the single edge crack problem: (a) 603 elements, (b) 697 elements, and (c) 859 elements.}
	\label{Fig.11m}
\end{figure}

\begin{figure}[!httb]
	\centering
	\subfigure[]{
		\includegraphics[scale=1]{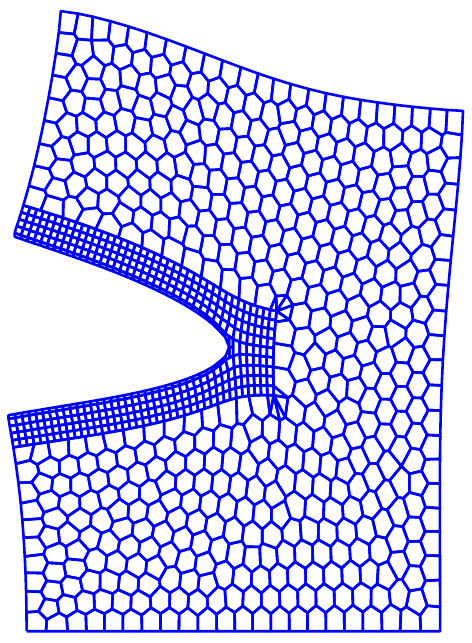}
		\label{Fig.12a}
	} \hspace*{4em}
	\subfigure[]{
		\includegraphics[scale=1]{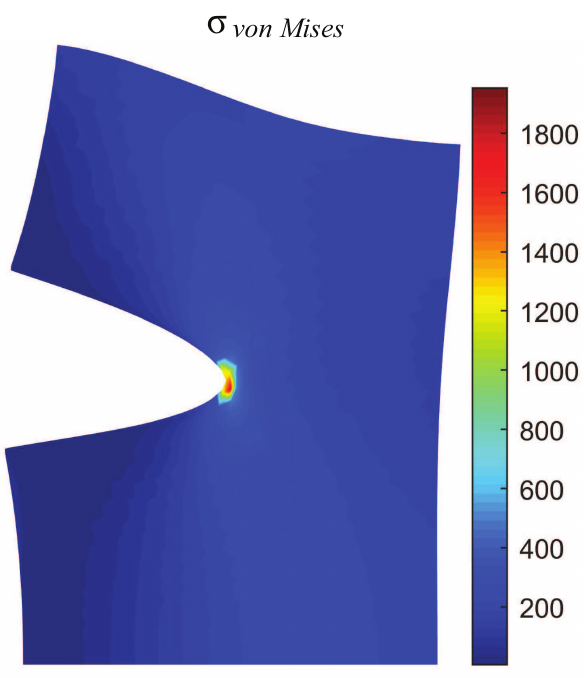}
		\label{Fig.12b}
	}
	\caption{(a) Mesh deformation, (b) distributions of von Mises stress.}
	\label{Fig.12}
\end{figure}

Consider a square domain of size $L = 2$ mm with an edge crack of length $a = 1$ mm as shown in Fig. \ref{Fig.11}. The bottom edge is fixed in y-direction, and the node at the left bottom corner is blocked to prevent a rigid body motion. A uniform load in the vertical direction is imposed on the edge. The domain is discretized into two kinds of meshes, including a polygonal mesh with three different local refinements (603, 697 and 859 elements) as shown in Fig. \ref{Fig.11m}, and a quadrilateral structured mesh with $49\times 49$ elements. A material parameter with $\mu = 0.4225$ MPa adopted to the Neo-Hookean incompressible material model in Eq. \ref{Eq_strain_ener_incom} is employed in the analysis. The computational process is executed in 40 load steps with each increment $\Delta\sigma = 5000$ Pa. The mesh deformation and stress field of von Mises are shown in Fig. \ref{Fig.12}, which exhibits concentrated stress in the vicinity of the crack tip.\\

\begin{figure}
	\centering
	\includegraphics[scale=1]{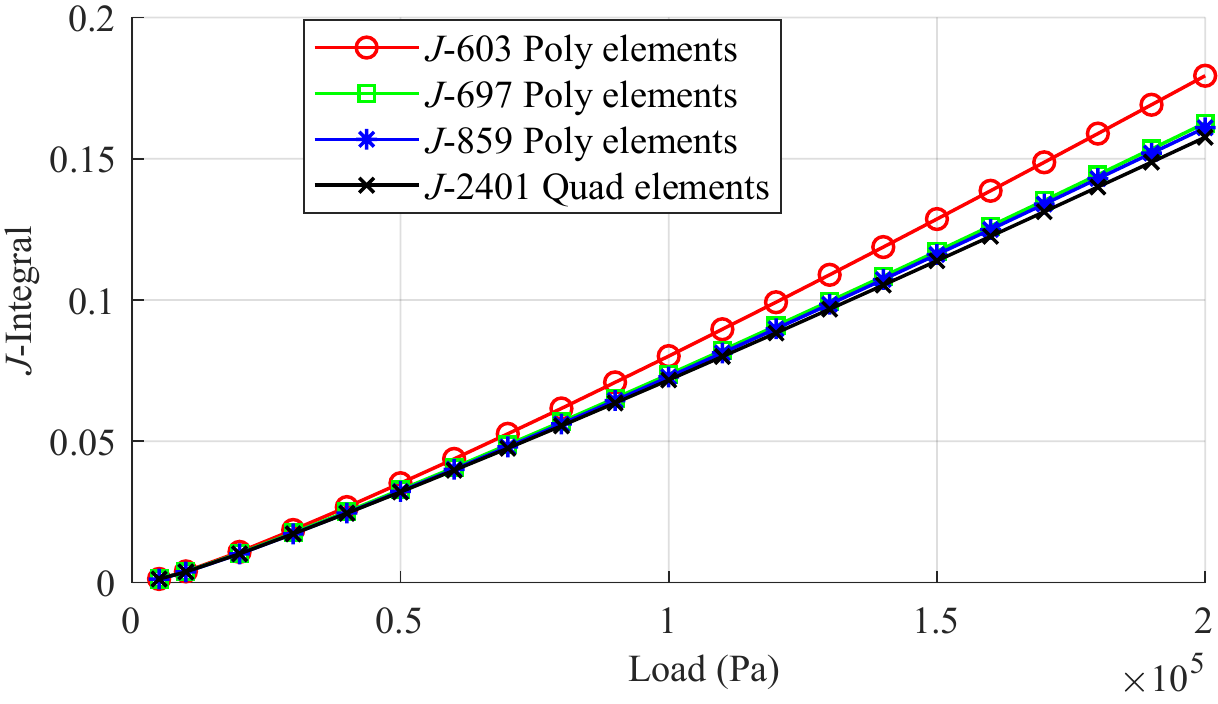}
	\caption{$J$-Integral obtained by the use of polygonal meshes and quadrilateral mesh.}
	\label{Fig.13}
\end{figure}

The values of the $J$-Integral obtained from the present method under the use of three meshes are plotted in Fig. \ref{Fig.13}. It should be noted that the evaluation of the $J$-Integral is stable when the integration domain radius is chosen as $r = 3*{h_{size}}$, see Refs. \citep{Belytschko_Black_1999,Legrain_Moes_Verron_2005,Rashetnia_2015}. A comparison of $J$ values analyzed from the three polygonal meshes and quadrilateral one shows that the performance of good refinements is lower than that of the coarse refinement. In particular, the figure obtained from both meshes (697 and 859 elements) is different with that of quadrilateral mesh with only $5\%$ while that from 603-element mesh is over $10 \%$.\\

\begin{figure}
	\centering
	\includegraphics[scale=1]{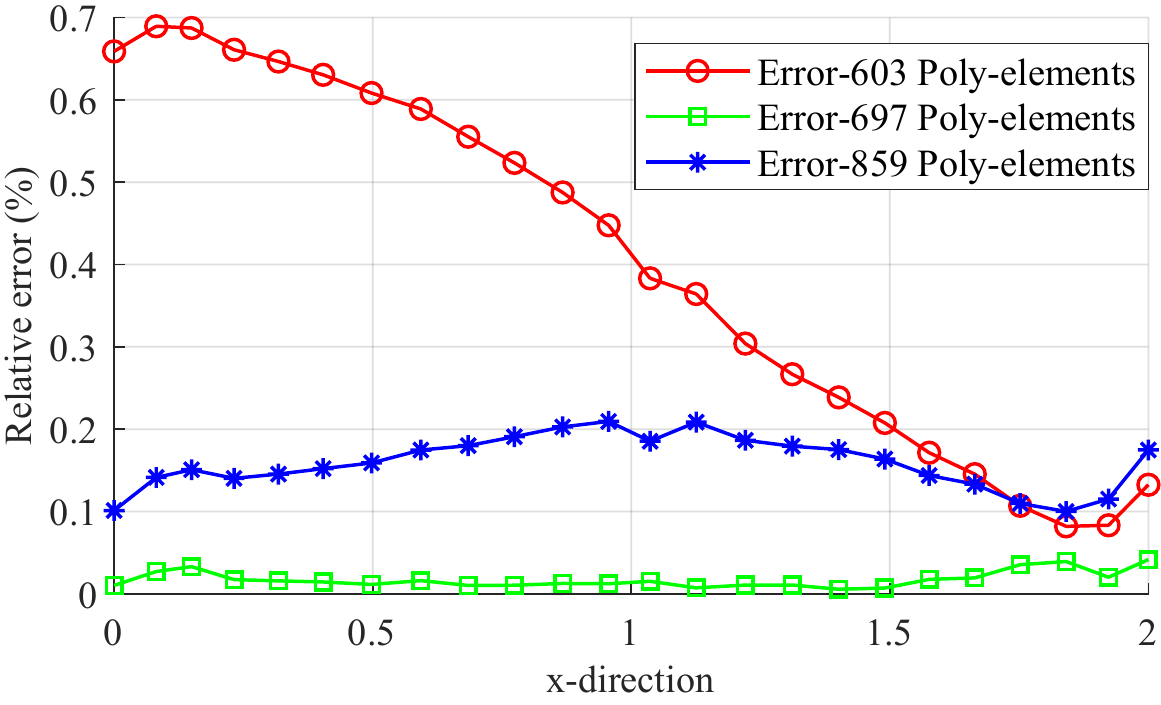}
	\caption{Relative error of the displacements along top edge.}
	\label{Fig.13ad}
\end{figure}

Fig. \ref{Fig.13ad} shows the relative error of the displacement along the top edge regarding the reference solution obtained from the use of quadrilateral mesh. The performance from the polygonal mesh with a very fine refinement is much lower than the counterpart. From displays in Fig. \ref{Fig.13} and Fig. \ref{Fig.13ad}, the local refinement has significant effects on the computational results. With only 697 elements, the obtained results well fit those from a very fine structured mesh with 2041 elements.

\subsection{A center crack in a square plate}
\label{subsec 4.2 center crack}

\begin{figure}[!httb]
	\centering
	\subfigure[]{
		\includegraphics[scale=1]{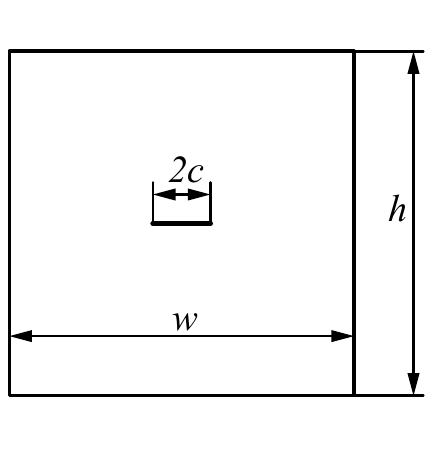}
		\label{Fig.center_geo}
	} \hspace*{0em}
	\subfigure[]{
		\includegraphics[scale=1]{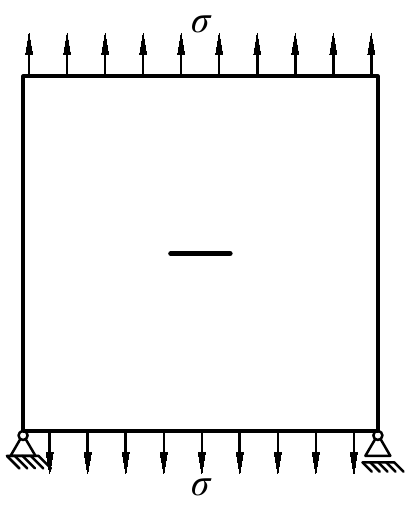}
		\label{Fig.center_axial}
	} \hspace*{0em}
	\subfigure[]{
		\includegraphics[scale=1]{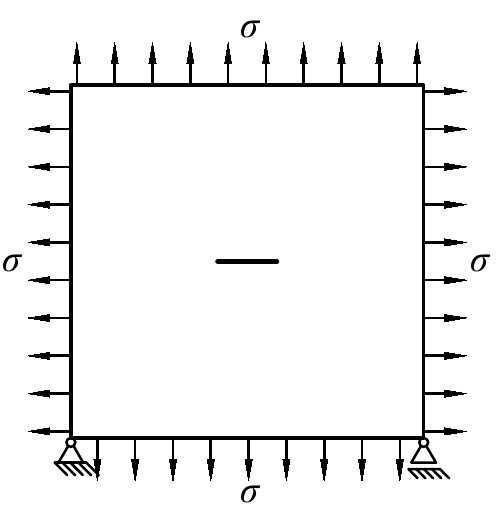}
		\label{Fig.center_biaxial}
	}
	\caption{A center crack problem: (a) a center crack domain, (b) uniaxial extension, (c) equibiaxial extension.}
	\label{Fig.centercrack}
\end{figure}

\begin{figure}[!httb]
	\centering
	\subfigure[]{
		\includegraphics[scale=1]{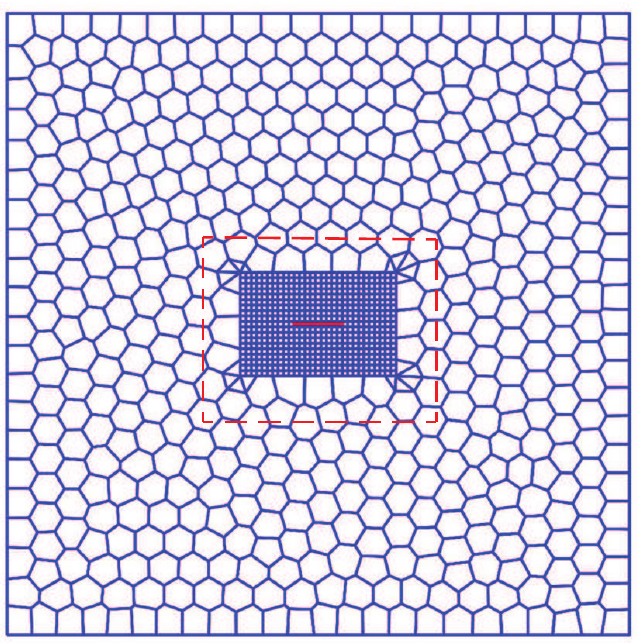}
		\label{Fig.mesh_center}
	} \hspace*{0em}
	\subfigure[]{
		\includegraphics[scale=1]{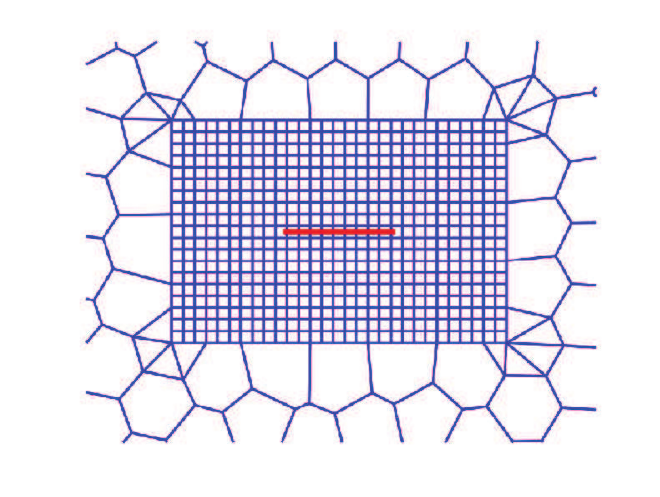}
		\label{Fig.mesh_zoom}
	}
	\caption{(a) Polygonal mesh (1066 elements), (b) a zoom around the crack path.}
	\label{Fig.Polymesh_center}
\end{figure}

\begin{figure}
	\centering
	\includegraphics[scale=1]{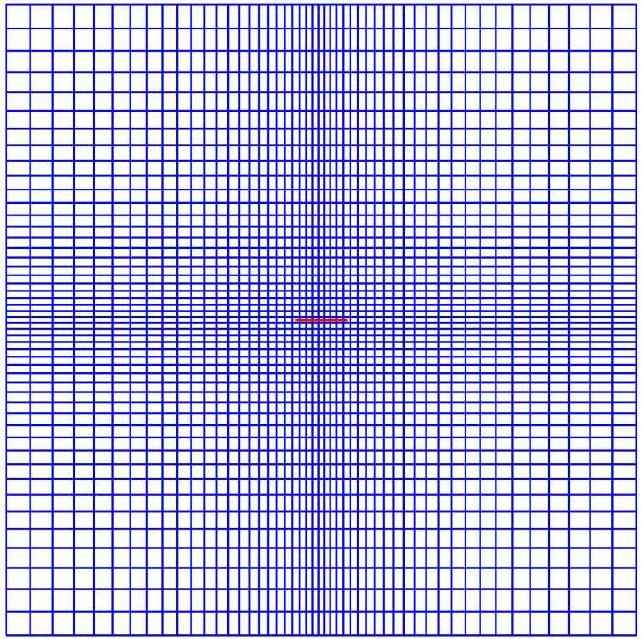}
	\caption{Quadrilateral mesh (2401 elements).}
	\label{Fig.center_Q4mesh}
\end{figure}

The example is examined according to the Griffith problem to evaluate behaviors of the nonlinear elastic fracture mechanics. Fig. \ref{Fig.center_geo} illustrates the size of the domain given $h = w = 6$ mm with the presence of a crack of length $2c = 0.5$ mm. The benchmark is imposed by two kinds of uniform loading, including uniaxial extension and equibiaxial extension as shown in Fig. \ref{Fig.center_axial} and Fig. \ref{Fig.center_biaxial}, respectively, and the boundary condition for preventing rigid body motions is also supported for both. The computational domain is discretized by a modified polygonal mesh with 1066 elements which consists of a very fined local quadrilateral mesh shown in  Fig. \ref{Fig.Polymesh_center}. A quadrilateral mesh with a high element density around the crack path as depicted in Fig. \ref{Fig.center_Q4mesh} is also considered in the analysis. The material is chosen to be incompressible with shear modulus $\mu  = 0.4225$ MPa. The main investigation of this problem is employed in the evolution of the $J$-Integral corresponding to each stretch level. The obtained results of the $J$-integral are compared with values of the tearing energy introduced by Lake \citep{Lake_1970}, and Lindley \citep{Lindley_1972} for the center crack problem in axial extension and Yeoh \citep{Yeoh_2002} in equibiaxial extension. In detail, an approximation factor $k$ dependent on the principal extension $\lambda$ is proposed to assess the tearing energy. The following form is given by
\begin{equation}
          G = 2kWc
          \label{Eq.tearing_ener}
\end{equation}
where $W$ is the strain energy density of the computational domain without the crack.\\

From Lake's work, the factor $k$ is measured as
\begin{equation}
         k = \dfrac{\pi }{{\sqrt \lambda  }}
         \label{Eq.Lake}
\end{equation}

and from the assumption of Lindley, the factor $k$ is
\begin{equation}
         k = \frac{{2.95 - 0.08\left( {1 - \lambda } \right)}}{{\sqrt \lambda  }}
         \label{Eq.Lindley}
\end{equation}

\begin{figure}
	\centering
	\includegraphics[scale=1]{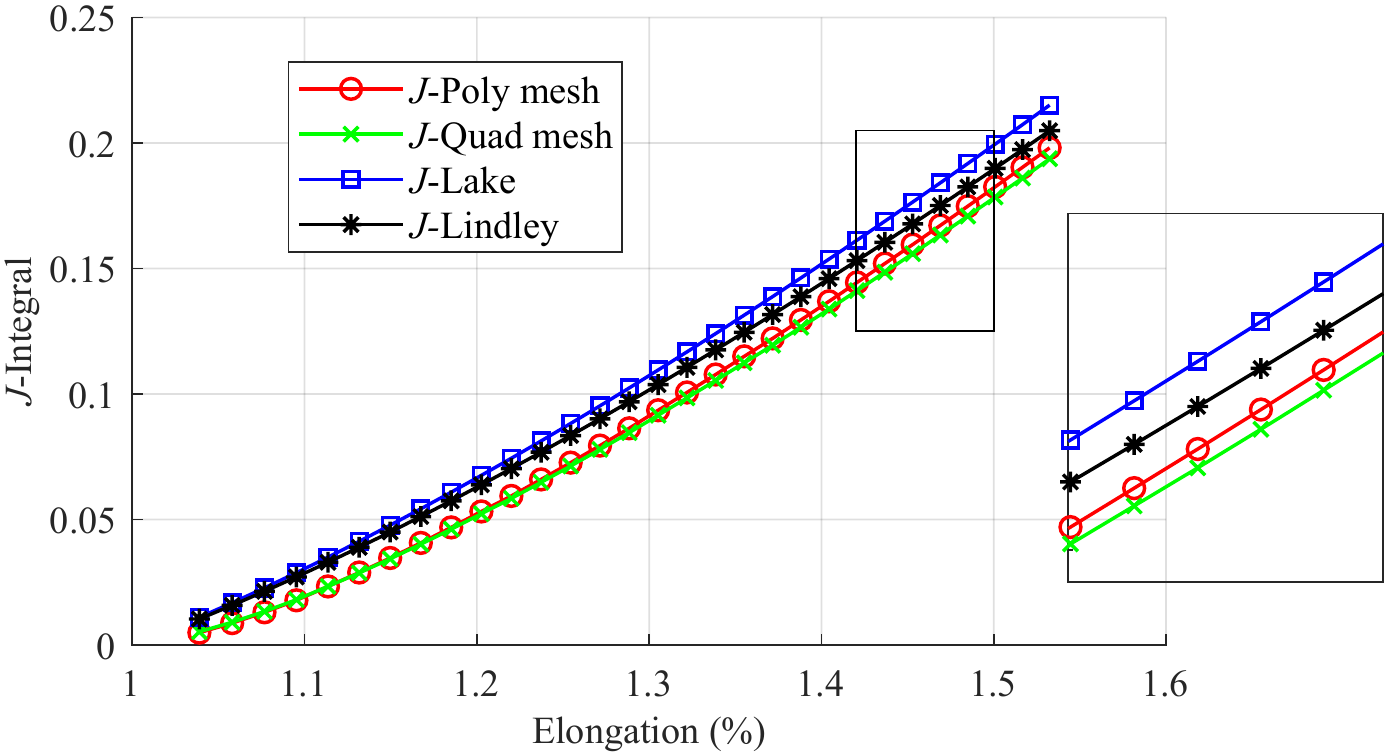}
	\caption{The variation of $J$-Integral in the case of uniaxial extension.}
	\label{Fig.J_axial_center}
\end{figure}

Fig. \ref{Fig.J_axial_center} shows the comparison of the $J$-Integral’s evolution with respect to elongation levels obtained from the present numerical method and two analytical approaches. As expected, the numerical results become close to the Lindley’s curve and well fit at higher stretch levels. In a close-up view of this figure, the performance of polygonal mesh shows a higher accuracy than that of quadrilateral mesh.\\

\begin{figure}
	\centering
	\includegraphics[scale=1]{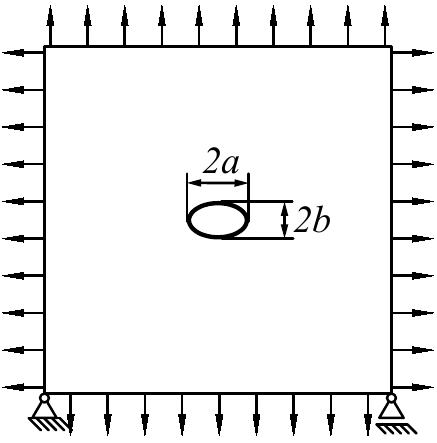}
	\caption{Assumption of the crack's deformation.}
	\label{Fig.crack_deform}
\end{figure}

The center crack problem is successively developed in the case of equibiaxial extension by assuming a crack with the length $c$ to become an ellipse as depicted in Fig. \ref{Fig.crack_deform}. The numerical factor is given by
\begin{equation}
          k = \frac{{{\sigma _y}\pi b}}{{4Wc}}
          \label{Eq.J_Yeoh}
\end{equation}
in which the stress $\sigma_y$ as reported in \citep{Yeoh_2002} is computed as
\begin{equation}
          \sigma _y = 2\left( {{\lambda _y} - \lambda _y^{ - 5}} \right)\left( {\frac{{\partial W }}{{\partial {I_1}}} + \lambda _y^2\frac{{\partial W }}{{\partial {I_2}}}} \right)
          \label{Eq.sigma_Yeoh}
\end{equation}
where $I_1$ and $I_2$ are the principal invariants.

\begin{figure}
	\centering
	\includegraphics[scale=1]{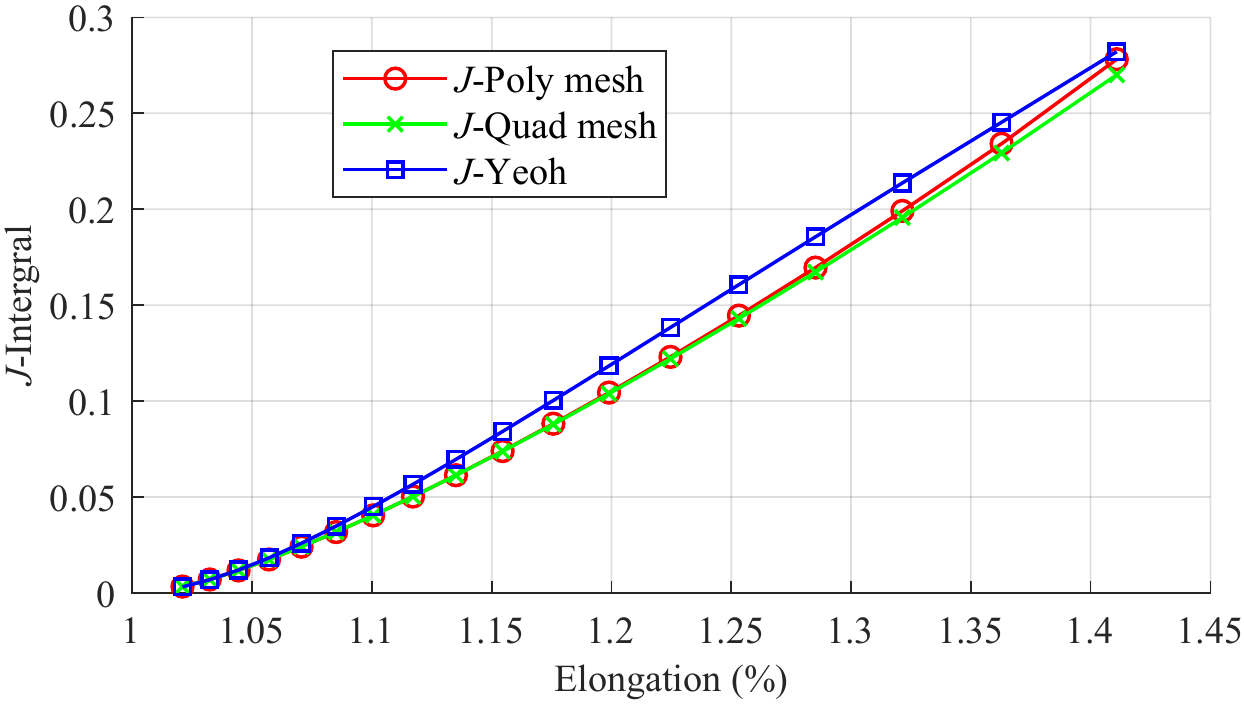}
	\caption{The variation of $J$-Integral in the case of equibiaxial extension.}
	\label{Fig.J_biaxial_center}
\end{figure}

The $J$-Integral values obtained from the present method with respect to variation of the axial extension are plotted in Fig. \ref{Fig.J_biaxial_center} and are also compared with those from the Yeoh's approach. We observe that when the axial stretch becomes large, the numerical $J$-Integral's performance tends to fit the result computed by the Yeoh's path.
\\

To address capacity of the present approach being generalized to any hyper-elastic models, we additionally consider other hyperelastic model named as the Mooney-Rivlin model. This model can be seen as a special case of the well-known Ogden model in which the strain energy for incompressible materials can be expressed in the following general form as
\begin{equation}
	W \left( {\bf{C}} \right) = \sum\limits_{p = 1}^N {\frac{{{\mu _p}}}{{{\alpha _p}}}} \left( {\lambda _1^{{\alpha _p}} + \lambda _1^{{\alpha _p}} + \lambda _1^{{\alpha _p}} - 3} \right)
	\label{Eq.Ogden}
\end{equation}
where $\lambda _1$, $\lambda _2$, $\lambda _3$ are principle stretches, and $(\mu _p,\, \alpha _p)$ are material parameters.\\

The Mooney-Rivlin model can be obtained by setting $N=2$, ${\alpha _1} = 2$, and ${\alpha _2} = -2$ from Eq. \ref{Eq.Ogden}
\begin{equation}
	W\left({\bf{C}} \right) = \frac{{{\mu _1}}}{2}\left( {{I_1} - 3} \right) - \frac{{{\mu _2}}}{2}\left( {{I_2} - 3} \right)
	\label{Eq.Mooney}
	\end{equation}
%	where ${I_1} = trace({\bf{C}})$, and ${I_2} = \dfrac{1}{2} \left( {{{\left( {trace({\bf{C}})} \right)}^2} - trace({{\bf{C}}^2})} \right)$.
\\

\begin{figure}[!httb]
	\centering
	\subfigure[]{
		\includegraphics[scale=1]{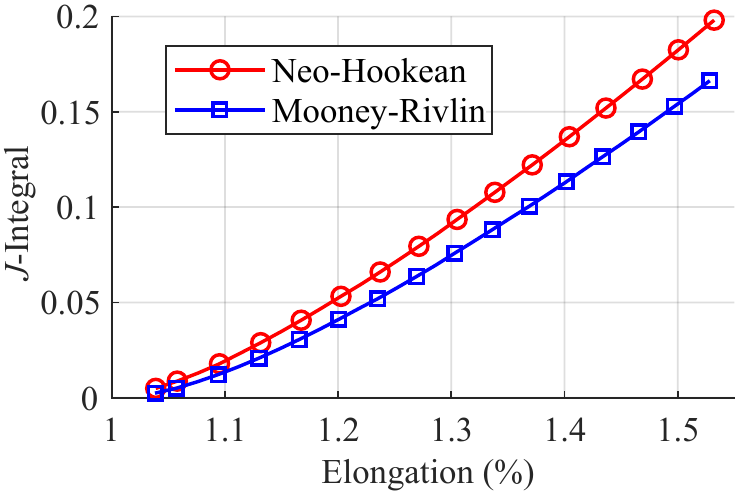}
		\label{Fig.Mooney_Neo_axial}
	} %\hspace*{1em}
	\subfigure[]{
		\includegraphics[scale=1]{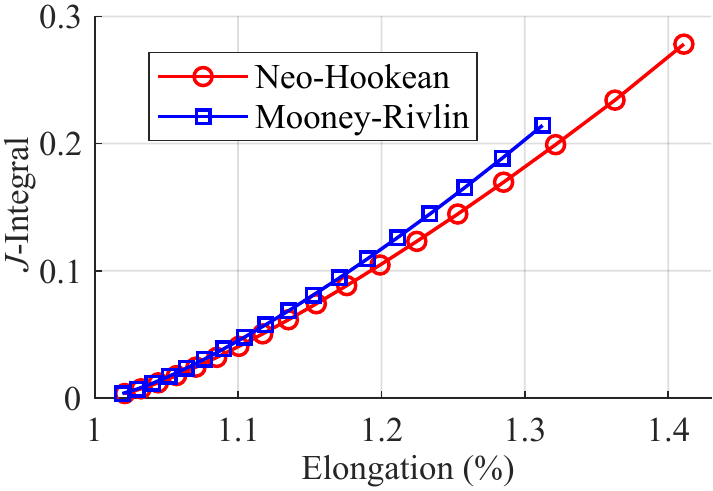}
		\label{Fig.Mooney_Neo_biaxial}
	}
	\caption{Comparison of $J$-Integral between two hyperelastic models: (a) uniaxial extension, (b) equibiaxial extension.}
	\label{Fig.J_Mooney_Neo}
\end{figure}

For computation, material parameters are chosen consistently to the input data of Neo-Hookean incompressible material given. In particular, they are typically given by ${\mu _1} = 3.6969 \times {10^5}$ Pa, and ${\mu _2} = -0.5281 \times {10^5}$ Pa. The numerical results are obtained for both cases of uniaxial and equibiaxial extension, and loading magnitudes are similar to the former counterparts. The assessment of $J$-Integral for Mooney-Rivlin model is shown in Fig. \ref{Fig.J_Mooney_Neo} and is compared with that of the Neo-Hookean model. It is seen that the release energy attained from the Mooney-Rivlin model for both stretching cases is different from the Neo-Hookean one although the material
is managed to be a typical rubber material. As depicted in Fig. \ref{Fig.J_Mooney_Neo}, the present approach is well capable of evaluating some hyperelastic models when the corresponding enrichment function to each hyperelastic material model is defined.
\subsection{A single edge crack specimen in a rectangular plate}
\label{sub_sec 4.2}

\begin{figure}[!httb]
	\centering
	\subfigure[]{
		\includegraphics[scale=1]{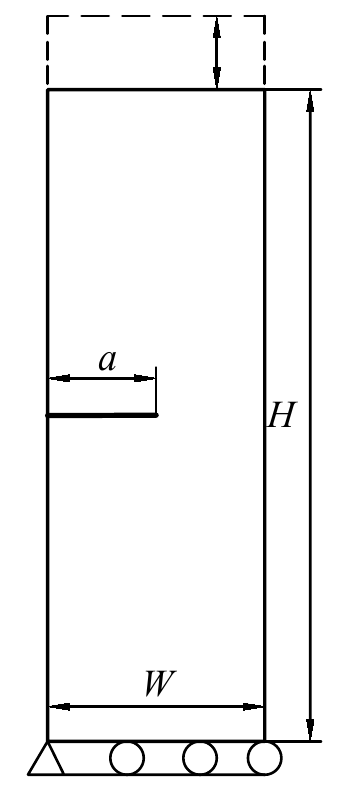}
		\label{Fig.domain_edge_rec}
	} \hspace*{5em}
	\subfigure[]{
		\includegraphics[scale=1]{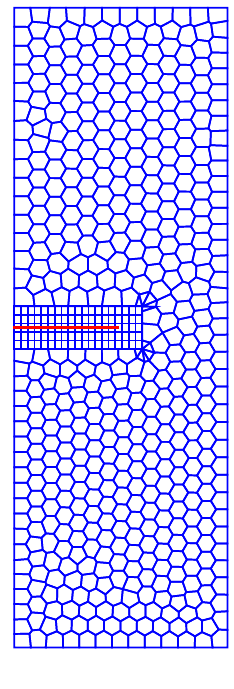}
		\label{Fig.mesh_edge_rec}
	}
	\caption{A single edge-crack problem in a rectangle plate: (a) domain geometry and boundary conditions, (b) polygonal mesh.}
	\label{Fig.14}
\end{figure}

A further single edge crack problem is continuously taken into the consideration of the effect of geometrically nonlinear structures. The boundary conditions at the bottom edge are subjected to the same manner in the previous example \ref{sub_sec 4.1}. The top edge is pulled in tension via a controlled displacement as shown in Fig. \ref{Fig.domain_edge_rec}. The size of the domain is chosen to be $W = 2$ mm, and $H = 6$ mm, with a crack length of $1$ mm. A mesh with 604 elements is employed for the analysis as shown in Fig. \ref{Fig.mesh_edge_rec}. The material is chosen to be nearly incompressible, and is analyzed based on the strain energy form in Eq. \ref{Eq_strain_ener_compress}. The material's parameters are given as $E = 50\times 10^3\,N/{mm^2}$, and $\nu = 0.45$, and the plane train is assumed. The process is executed in 20 steps until the last elongation of $16.6\%$ is obtained.\\

\begin{figure}[!httb]
	\centering
	\subfigure[]{
		\includegraphics[scale=1]{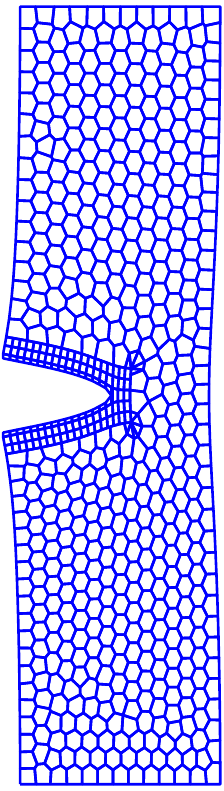}
		\label{Fig.edge_rec_deform_non}
	} \hspace*{1em}
	\subfigure[]{
		\includegraphics[scale=1]{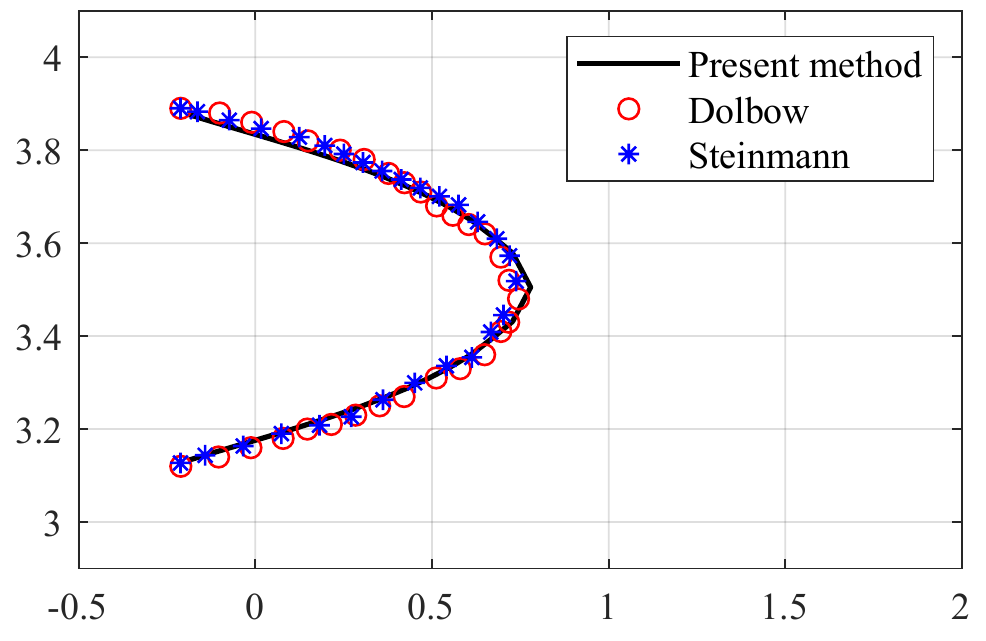}
		\label{Fig.curve_crack_rec_non}
	}
	\caption{(a) Mesh deformation for nonlinear fracture, (b) comparison of the crack with published results.}
	\label{Fig.edge_crack_non}
\end{figure}

This benchmark is implemented in both cases of linear and nonlinear fracture problem. An explicit result regarding the mesh finite-deformation is shown in Fig. \ref{Fig.edge_rec_deform_non}. The deformation of the crack path as illustrated in Fig. \ref{Fig.curve_crack_rec_non} is separately plotted and compared with the previous results found in \citep{Dolbow_Devan_2004,Steinmann_Ackermann_Barth_2001}.\\

\begin{figure}[!httb]
	\centering
	\subfigure[]{
		\includegraphics[scale=1]{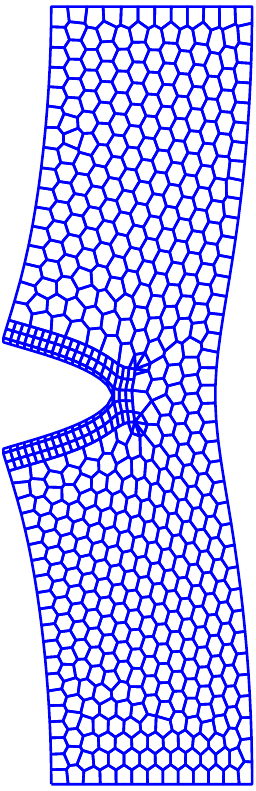}
		\label{Fig.edge_rec_deform_lin}
	} \hspace*{1em}
	\subfigure[]{
		\includegraphics[scale=1]{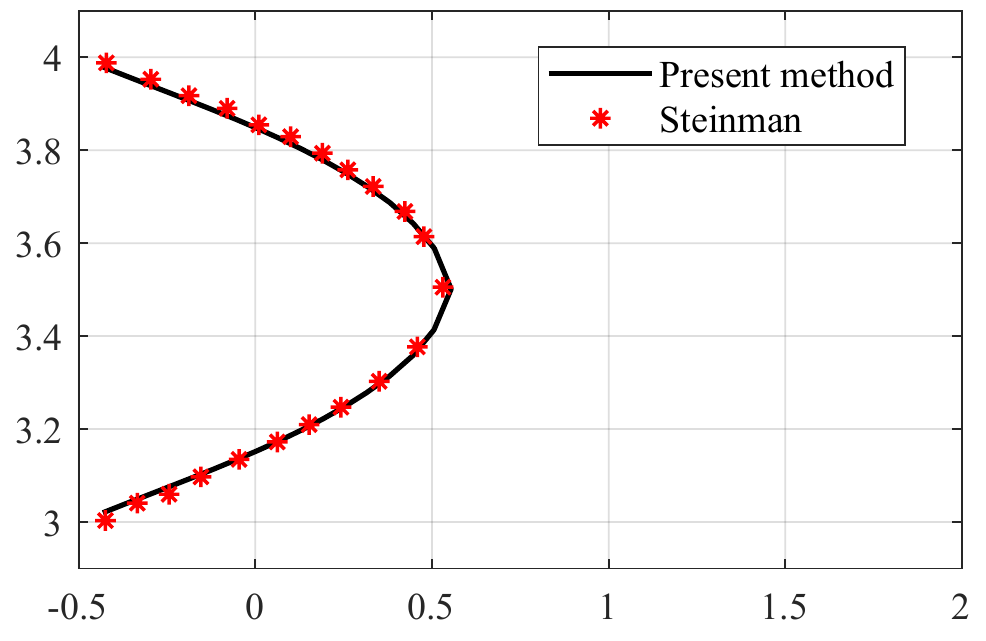}
		\label{Fig.curve_crack_rec_lin}
	}
	\caption{(a) Mesh deformation for linear fracture, (b) comparison of the crack with published result.}
	\label{Fig.edge_crack_lin}
\end{figure}

Fig. \ref{Fig.edge_crack_lin} indicates the deformed configuration of the domain  under the linear fracture analysis. The performance acquired from the present method is favorably compared with that achieved by Steinman, reported in \citep{Steinmann_Ackermann_Barth_2001}. Obviously, the obtained results for both cases in fracture mechanics are in good agreement with the published ones.
%In addition, a smooth visualization for the stress distributions of vonMises shown in Fig. \ref{Fig.15c} demonstrates the advantages of the polygonal finite elements for evaluating stresses.

%Another investigation is the numerical computation of the $J$-integral, whose SIFs values are normalized to make a comparison with those of the linear fracture problem. With the imposed boundary conditions, the SIFs obtained from the $J$-integral are evaluated only for the pure mode $I$. The results of the $J$-integral in the nonlinear analysis and those in the linear counterpart are numerically calculated with the same extension levels. A comparison of the $J$-integral obtained the two frameworks is shown in Fig. \ref{Fig.16}. The performance of the nonlinear analysis appears in a lower value compared with that of counterpart. This demonstrates that the release energy represented as the $J$-integral is lower than the linear achievement.
%\begin{figure}
%	\centering
%	\includegraphics[scale=1]{F16_ben2_SIF.pdf}
%	\caption{A comparison of J-integral between the linear and nonlinear analysis.}
%	\label{Fig.16}
%\end{figure}

\subsection{An inclined edge crack in a finite plate}
\label{sub_sec 4.3}
\begin{figure}
	\centering
	\includegraphics[scale=1]{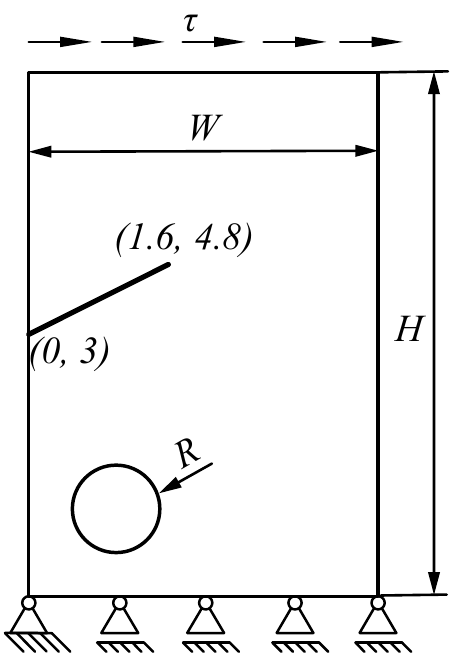}
	\caption{A inclined edge crack problem.}
	\label{Fig.inclined_problem}
\end{figure}

\begin{figure}[!httb]
	\centering
	\subfigure[]{
		\includegraphics[scale=1.01]{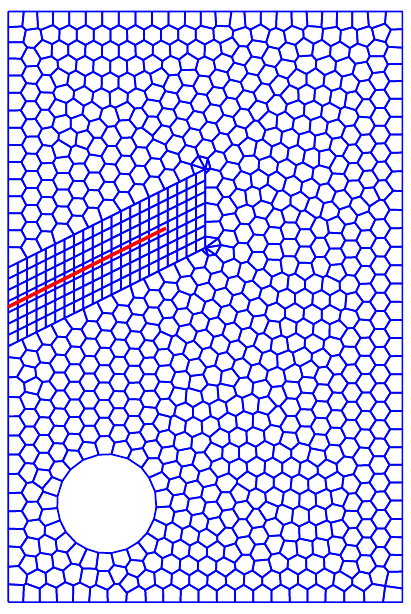}
		\label{Fig.17a}
	} \hspace*{5em}
	\subfigure[]{
		\includegraphics[scale=1]{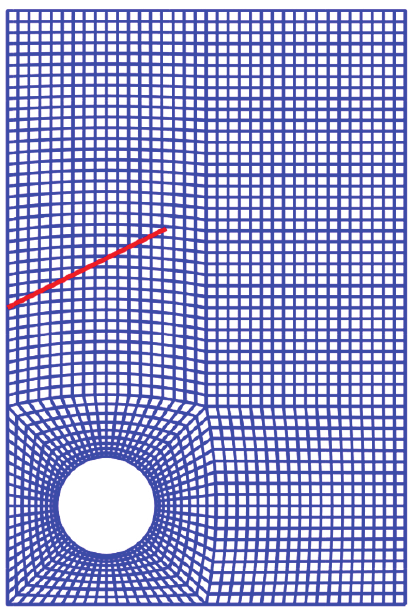}
		\label{Fig.17b}
	}
	\caption{Mesh discretization: (a) Polygonal mesh, (b) quadrilateral mesh.}
	\label{Fig.17}
\end{figure}

This example investigates the geometrical changes and evolution of the $J$-Integral in a mixed-mode crack problem. An inclined edge crack in a rectangular plate with a circular hole is described in Fig. \ref{Fig.inclined_problem}, in which their definitions are given by $W=4$ mm, $H=6$ mm, and $R=0.5$ mm. A shear loading $\tau$ is imposed to the top edge, and prescribed boundary conditions are installed at the bottom edge. A compressible material is chosen to be in the plane strain condition with Young's modulus $E = 50\times 10^3\,N/{mm^2}$ and Poisson's ratio $\nu = 0.3$. A polygonal mesh of 955 elements and a fine quadrilateral mesh of 2268 elements as depicted in Fig. \ref{Fig.17} are employed in analysis. The computational process is carried out with 16 load steps under each increment $\Delta\tau=50$ N.\\

\begin{figure}
	\centering
	\includegraphics[scale=1]{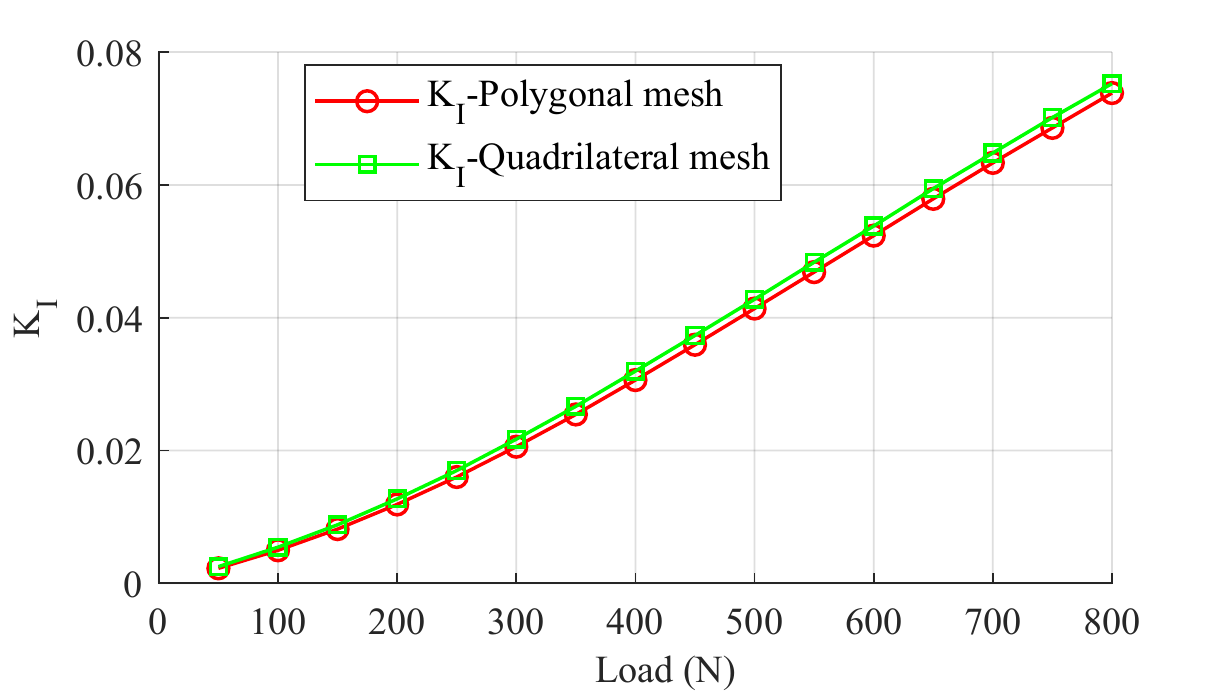}
	\caption{Stress intensity factors of mode $I$.}
	\label{Fig.18}
\end{figure}

\begin{figure}
	\centering
	\includegraphics[scale=1]{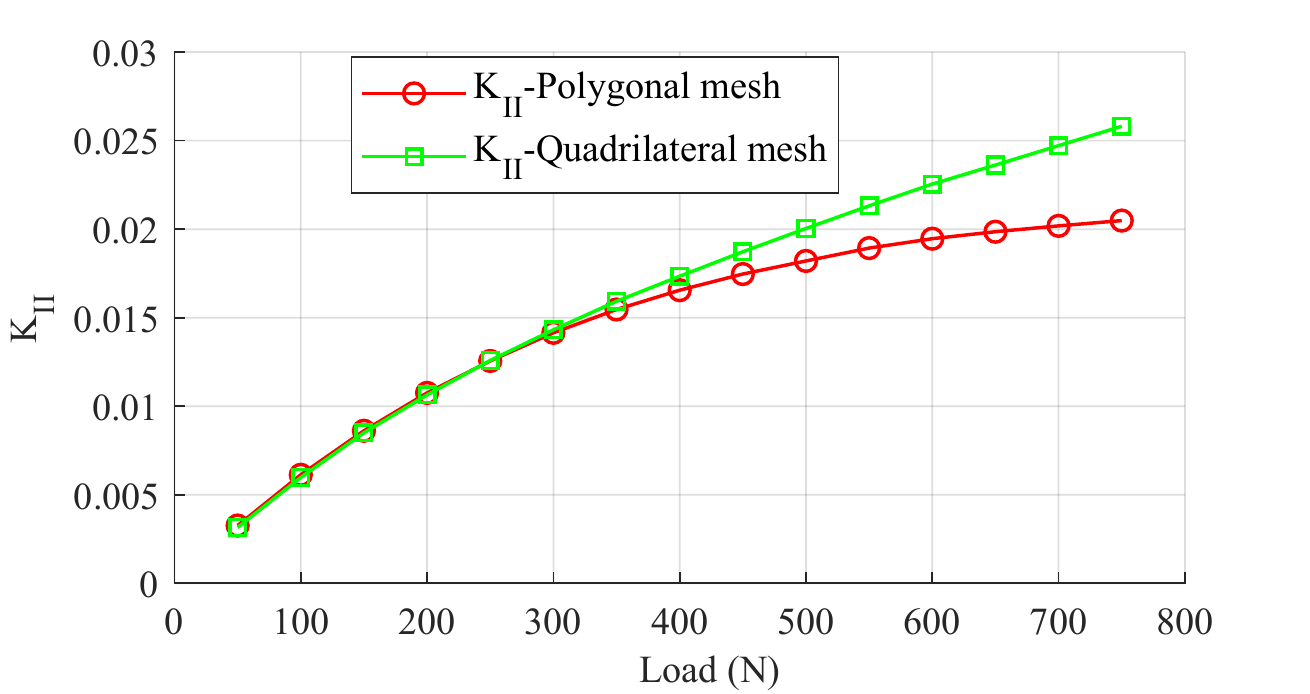}
	\caption{Stress intensity factors of mode $II$.}
	\label{Fig.19}
\end{figure}

A subsequent numerical examination is taken into the evaluation of the $J$-Integral to get values of the SIFs for each mode. The results obtained under the use of two finite mesh approaches including polygonal and quadrilateral meshes are favorable compared to tackle the computational efficiency of the present XFEM framework. In particular, the evolution of SIFs of mode $I$ is shown in Fig. \ref{Fig.18} which explicitly exhibits a very good agreement of the performance from the polygonal mesh with the quadrilateral counterpart while the number of polygonal elements is around a half of quadrilateral ones. Fig. \ref{Fig.19} indicates the comparison of the SIFs for mode $II$. It is seen that two approaches produce almost same results in first steps, and there is a difference in the rest which may be acceptable when it comes to the consideration of mode $II$.\\

\begin{figure}[!httb]
	\centering
	\subfigure[]{
		\includegraphics[scale=1]{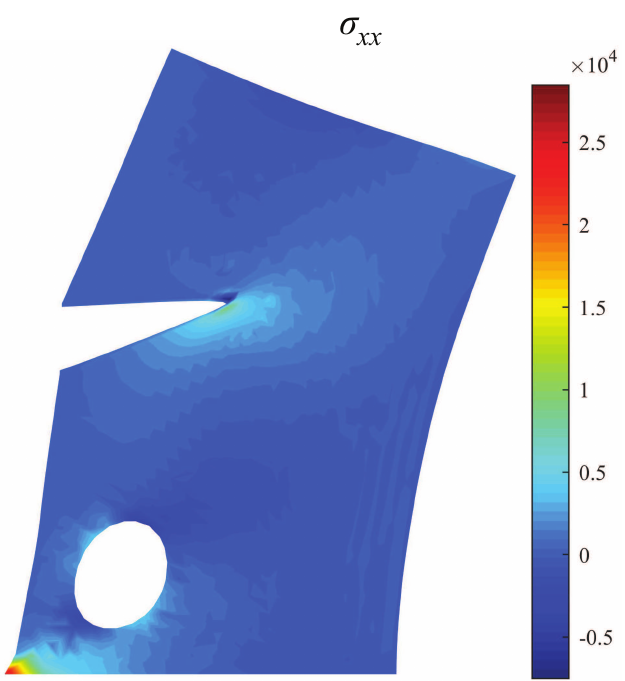}
		\label{Fig.sigmaxx}
	} \hspace*{2em}
	\subfigure[]{
		\includegraphics[scale=1]{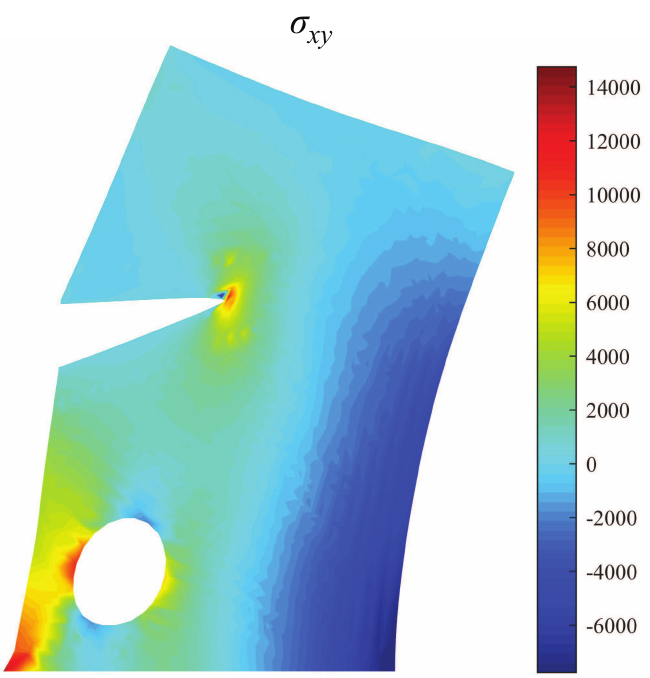}
		\label{Fig.sigmaxy}
	}
	\subfigure[]{
		\includegraphics[scale=1]{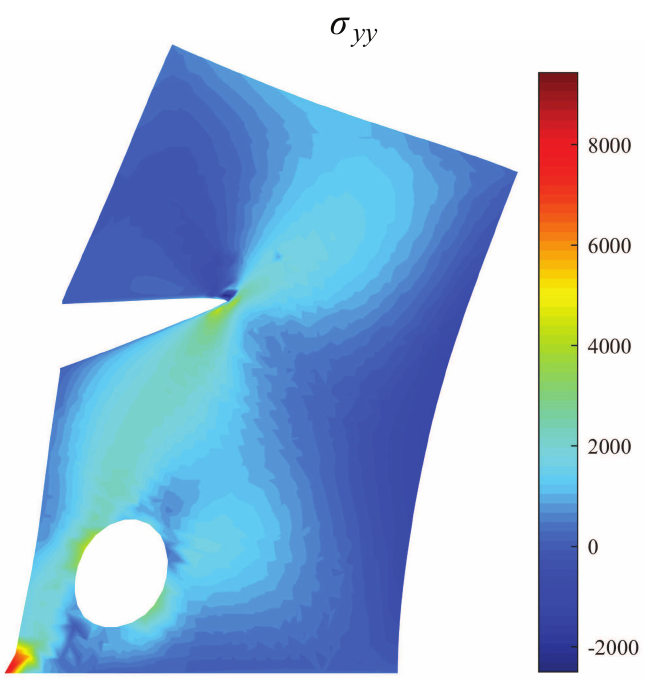}
		\label{Fig.sigmayy}
	} \hspace*{2em}
	\subfigure[]{
		\includegraphics[scale=1]{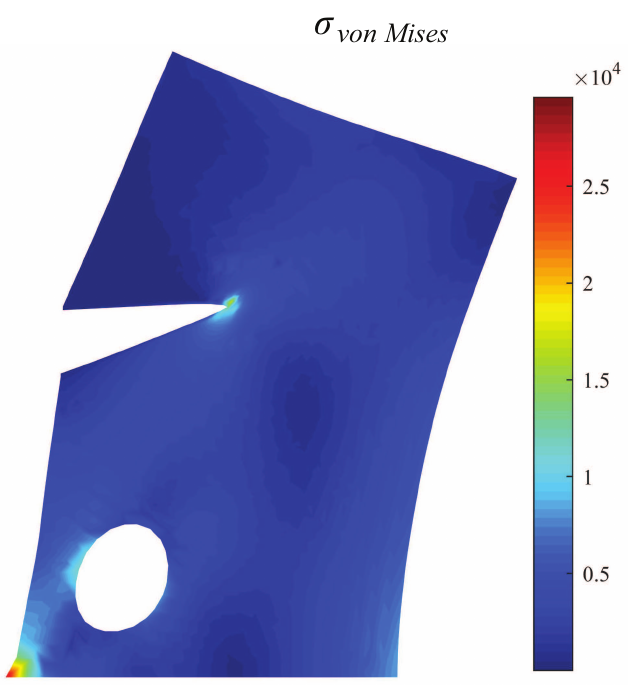}
		\label{Fig.sigmavm}
	}
	\caption{Stress distributions: (a) $\sigma_{xx}$, (b) $\sigma_{xy}$, (c) $\sigma_{yy}$, (b) $\sigma_{vonMises}$.}
	\label{Fig.20}
\end{figure}

Fig. \ref{Fig.20} shows stress distributions which are plotted in the deformed configuration of the crack domain. These figures provide a smooth visualization beneficial to observe the changes in the geometry of crack path, and hole, as well as in the stress fields.\\

\subsection{A mechanism specimen with an edge crack}
\label{sub_sec 4.4}

\begin{figure}[!httb]
	\centering
	\subfigure[]{
		\includegraphics[scale=1]{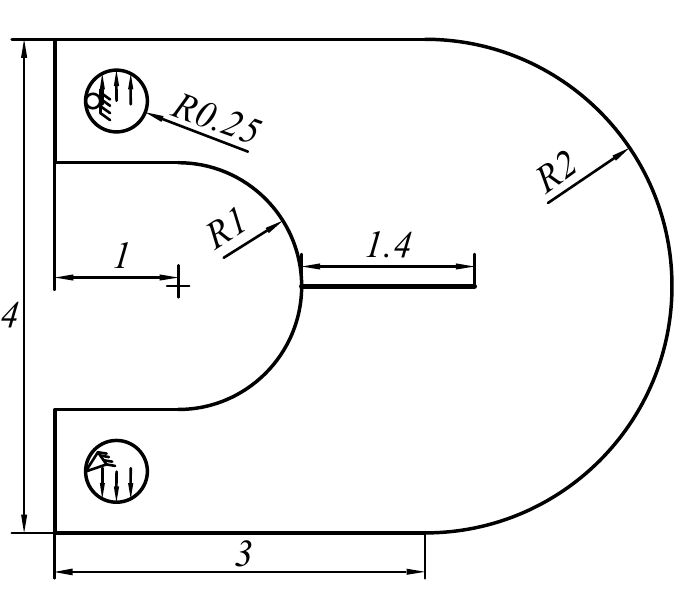}
		\label{Fig.21a}
	} \hspace*{0em}
	\subfigure[]{
		\includegraphics[scale=1]{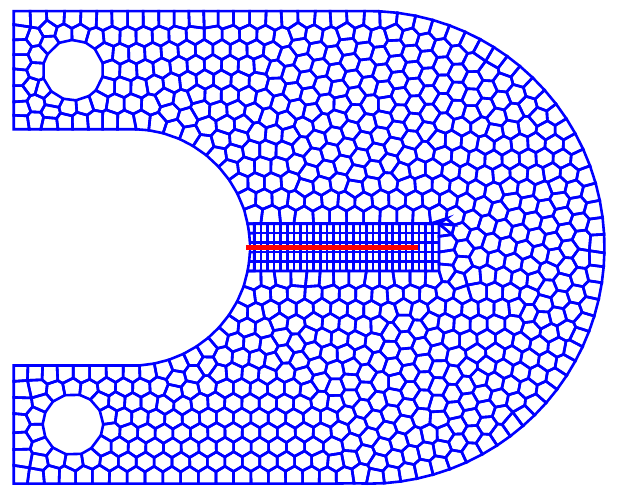}
		\label{Fig.21b}
	}
	\caption{(a) Mechanicsm with an edge crack, (b) modified mesh generation.}
	\label{Fig.21}
\end{figure}

This example examines a domain with a slightly complicated geometry input, which includes an edge crack. The problem geometry and boundary conditions are shown in Fig. \ref{Fig.21a}.  The specimen is subjected to a tensile stress $\sigma = 500$ N at the two holes, and material parameters are given similarly to the previous example \ref{sub_sec 4.3}. With the capacity of polygonal mesh generation, the problem geometry of the domain can be conformed well in the mesh of only 750 polygonal elements and a local refinement of 150 quadrilateral elements as shown in Fig. \ref{Fig.21b}.\\

\begin{figure}[!httb]
	\centering
	\subfigure[]{
		\includegraphics[scale=1]{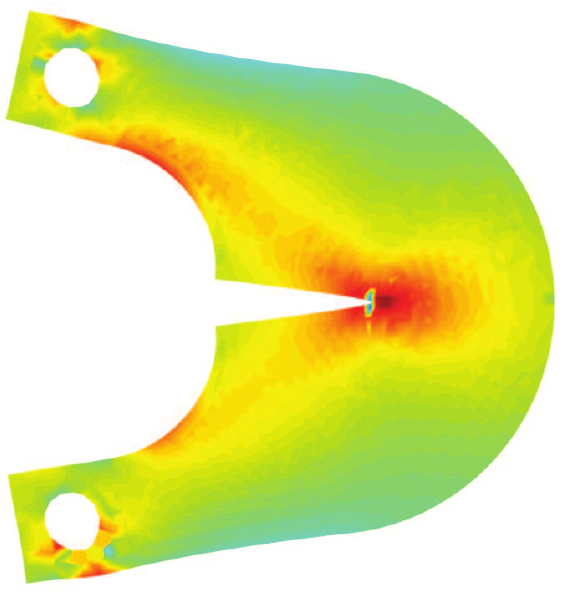}
		\label{Fig.22a}
	} \hspace*{0em}
	\subfigure[]{
		\includegraphics[scale=1]{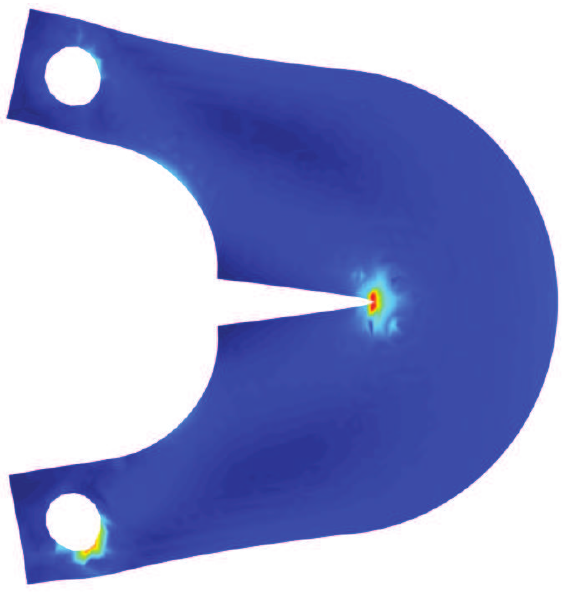}
		\label{Fig.22b}
	}
	\caption{Distributions: (a) the normal stress ${\sigma _{xx}}$, (b) the von Mises stress ${\sigma _{vonMises}}$.}
	\label{Fig.22}
\end{figure}

The performance of using polygonal elements to discretize the finite element space is sufficient for modeling complicated geometries with only a certain number of elements. It is because of using a fully unstructured polygonal mesh based on the algorithm of Voronoi cells, which provide a great flexibility in discretizing spaces. Moreover, the application of basis function over these elements owns a high-order character; thus calculations related to their derivatives are better than their counterpart with the of bilinear functions. The exhibition of stress distributions as shown in Fig. \ref{Fig.22} is a very smooth visualization for both the normal stress in x-direction and a the equivalent stress of von Mises.\\

\section{Conclusions}
\label{sec 5}

This study presented an extended polygonal finite element method to model the finite hyper-elastic deformation in fracture mechanics problems. The finite formulations were defined in the updated Lagrangian description appropriate to the numerical investigation of large strain problems. The use of a Ramp function was provided to handle the blending problems to enhance the accuracy of the XFEM technique. An introduction of the asymptotic enrichment function in the undeformed configuration and then a transformation from the undeformed to the deformed configuration were exhibited to capture the crack zone in the analysis process. A suitable enrichment function near the crack tip is taken into the fracture consideration of Neo-Hookean hyper-elastic material model. A novel approach of the space discretization on polygonal elements was offered to accordingly couple with local refinements on the structured quadrilateral mesh at the vicinity of the crack path. This improvement reduces the mesh distortions when a certain number of elements have significant changes in their geometry. With the capability of shape functions to conform to polygonal elements, the approximation fields over polygonal elements with possible hanging nodes are completely accomplished. In addition, there will be no limit to the number of hanging nodes additional located on original elements, which is very beneficial to the local refinement approach near the discontinuity areas.\\

Besides the noticeable points mentioned earlier, the need of further studies for the modeling of nonlinear fracture problems should be investigated. The asymptotic enrichment function in recent researches is only adapted to the Neo-Hookean material model, and thus developments for other nonlinear material models are necessary for the later examinations. The plasticity model is also supposed to make together with the present work in order to get general insights into the material behaviors and the nonlinear crack problems.

\section*{Acknowledgments}
This research is funded by Vietnam National Foundation for Science and Technology Development (NAFOSTED) under grant number 107.01-2018.32. The support provided by RISE-project BESTOFRAC (734370)–H2020 is gratefully acknowledged.

\appendix
\section{Definitions of Jacobian matrix}
\label{appendixA}

Components of the Jacobian matrix reported in Eq. \ref{Eq26} in the XFEM framework are defined as
\begin{small}
\begin{equation}
     \begin{array}{l}
     {J^{11}} = \dfrac{{\partial x}}{{\partial \xi }} = \sum\limits_{i \in I} {\dfrac{{\partial {N_i}}}{{\partial \xi }}{{\bar x}_i}}  + \sum\limits_{j \in J} {\dfrac{{\partial {N_j}}}{{\partial \xi }}\bar H({\bf{x}}){{{\bf{\bar a}}}_j}}  + \sum\limits_{k \in K} {\dfrac{{\partial {N_k}}}{{\partial \xi }}\sum\limits_l {{{\bar A}_l}({\bf{X}})R({\bf{X}})} {\bf{\bar b}}_k^l} \\
     + \sum\limits_{k \in K} {{N_k}\sum\limits_l {\left( {\dfrac{{\partial {{\bar A}_l}({\bf{X}})}}{{\partial X}}J_0^{11}R({\bf{X}}) + {{\bar A}_l}({\bf{X}})\dfrac{{\partial R({\bf{X}})}}{{\partial X}}J_0^{11} + \dfrac{{\partial {{\bar A}_l}({\bf{X}})}}{{\partial Y}}J_0^{12}R({\bf{X}}) + {{\bar A}_l}({\bf{X}})\dfrac{{\partial R({\bf{X}})}}{{\partial Y}}J_0^{12}} \right)} {\bf{\bar b}}_k^l}
     \end{array}
     \label{EqA1}
\end{equation}

\begin{equation}
     \begin{array}{l}
     {J^{12}} = \dfrac{{\partial y}}{{\partial \xi }} = \sum\limits_{i \in I} {\dfrac{{\partial {N_i}}}{{\partial \xi }}{{\bar y}_i}}  + \sum\limits_{j \in J} {\dfrac{{\partial {N_j}}}{{\partial \xi }}\bar H({\bf{x}}){{{\bf{\bar a}}}_j}}  + \sum\limits_{k \in K} {\dfrac{{\partial {N_k}}}{{\partial \xi }}\sum\limits_l {{{\bar A}_l}({\bf{X}})R({\bf{X}})} {\bf{\bar b}}_k^l} \\
     + \sum\limits_{k \in K} {{N_k}\sum\limits_l {\left( {\dfrac{{\partial {{\bar A}_l}({\bf{X}})}}{{\partial X}}J_0^{11}R({\bf{X}}) + {{\bar A}_l}({\bf{X}})\dfrac{{\partial R({\bf{X}})}}{{\partial X}}J_0^{11} + \dfrac{{\partial {{\bar A}_l}({\bf{X}})}}{{\partial Y}}J_0^{12}R({\bf{X}}) + {{\bar A}_l}({\bf{X}})\dfrac{{\partial R({\bf{X}})}}{{\partial Y}}J_0^{12}} \right)} {\bf{\bar b}}_k^l}
     \end{array}
     \label{EqA2}
\end{equation}

\begin{equation}
     \begin{array}{l}
     {J^{21}} = \dfrac{{\partial x}}{{\partial \eta }} = \sum\limits_{i \in I} {\dfrac{{\partial {N_i}}}{{\partial \eta }}{{\bar x}_i}}  + \sum\limits_{j \in J} {\dfrac{{\partial {N_j}}}{{\partial \eta }}\bar H({\bf{x}}){{{\bf{\bar a}}}_j}}  + \sum\limits_{k \in K} {\dfrac{{\partial {N_k}}}{{\partial \eta }}\sum\limits_l {{{\bar A}_l}({\bf{X}})R({\bf{X}})} {\bf{\bar b}}_k^l} \\
     + \sum\limits_{k \in K} {{N_k}\sum\limits_l {\left( {\dfrac{{\partial {{\bar A}_l}({\bf{X}})}}{{\partial X}}J_0^{21}R({\bf{X}}) + {{\bar A}_l}({\bf{X}})\dfrac{{\partial R({\bf{X}})}}{{\partial X}}J_0^{21} + \dfrac{{\partial {{\bar A}_l}({\bf{X}})}}{{\partial Y}}J_0^{22}R({\bf{X}}) + {{\bar A}_l}({\bf{X}})\dfrac{{\partial R({\bf{X}})}}{{\partial Y}}J_0^{22}} \right)} {\bf{\bar b}}_k^l}
     \end{array}
     \label{EqA3}
\end{equation}

\begin{equation}
     \begin{array}{l}
     {J^{22}} = \dfrac{{\partial y}}{{\partial \eta }} = \sum\limits_{i \in I} {\dfrac{{\partial {N_i}}}{{\partial \eta }}{{\bar y}_i}}  + \sum\limits_{j \in J} {\dfrac{{\partial {N_j}}}{{\partial \eta }}\bar H({\bf{x}}){{{\bf{\bar a}}}_j}}  + \sum\limits_{k \in K} {\dfrac{{\partial {N_k}}}{{\partial \eta }}\sum\limits_l {{{\bar A}_l}({\bf{X}})R({\bf{X}})} {\bf{\bar b}}_k^l} \\
     + \sum\limits_{k \in K} {{N_k}\sum\limits_l {\left( {\dfrac{{\partial {{\bar A}_l}({\bf{X}})}}{{\partial X}}J_0^{21}R({\bf{X}}) + {{\bar A}_l}({\bf{X}})\dfrac{{\partial R({\bf{X}})}}{{\partial X}}J_0^{21} + \dfrac{{\partial {{\bar A}_l}({\bf{X}})}}{{\partial Y}}J_0^{22}R({\bf{X}}) + {{\bar A}_l}({\bf{X}})\dfrac{{\partial R({\bf{X}})}}{{\partial Y}}J_0^{22}} \right)} {\bf{\bar b}}_k^l}
     \end{array}
     \label{EqA4}
\end{equation}
\end{small}

\section{Definitions of the matrix $\mathbf{B}$, and $\mathbf{G}$}
\label{appendixB}

As mentioned earlier, the matrix $\mathbf{B}$, and $\mathbf{G}$ in the XFEM are decomposed into 3 different parts, including standard part, possible enrichment parts for Heaviside and Asymptotic function. The definitions of these parts for node $I$ are given by

\begin{equation}
{\left( {{\bf{B}}_{}^{std}} \right)_I} = \left[ {\begin{array}{*{20}{c}}
	\vspace{1.5mm}{\dfrac{{\partial {N_I}}}{{\partial x}}}&0\\
	\vspace{1.5mm}0&{\dfrac{{\partial {N_I}}}{{\partial y}}}\\
	{\dfrac{{\partial {N_I}}}{{\partial y}}}&{\dfrac{{\partial {N_I}}}{{\partial x}}}
	\end{array}} \right]{\rm{ }}
\label{EqB1}
\end{equation}\\
\begin{equation}
     {\left( {{{\bf{B}}^H}} \right)_I} = \left[ {\begin{array}{*{20}{c}}
     	\vspace{1.5mm}{\dfrac{{\partial \left( {{N_I}\bar H({\bf{x}})} \right)}}{{\partial x}}}&0\\
     	\vspace{1.5mm}0&{\dfrac{{\partial \left( {{N_I}\bar H({\bf{x}})} \right)}}{{\partial y}}}\\
     	{\dfrac{{\partial \left( {{N_I}\bar H({\bf{x}})} \right)}}{{\partial y}}}&{\dfrac{{\partial \left( {{N_I}\bar H({\bf{x}})} \right)}}{{\partial x}}}
     	\end{array}} \right]{\rm{ }}
     \label{EqB2}
\end{equation}\\
\begin{equation}
     {\left( {{\bf{B}}_{}^{Tip}} \right)_I} = \left[ {\begin{array}{*{20}{c}}
     	\vspace{1.5mm}{\dfrac{{\partial \left( {{N_I}{{\bar A}_j}({\bf{X}})R({\bf{X}})} \right)}}{{\partial x}}}&0\\
     	\vspace{1.5mm}0&{\dfrac{{\partial \left( {{N_I}{{\bar A}_j}({\bf{X}})R({\bf{X}})} \right)}}{{\partial y}}}\\
     	{\dfrac{{\partial \left( {{N_I}{{\bar A}_j}({\bf{X}})R({\bf{X}})} \right)}}{{\partial y}}}&{\dfrac{{\partial \left( {{N_I}{{\bar A}_j}({\bf{X}})R({\bf{X}})} \right)}}{{\partial x}}}
     	\end{array}} \right]
     \label{EqB3}
\end{equation}

The term of $\mathbf{X}$ is the position defined in the initial configuration. The derivatives with respect to the current configuration are given by $\dfrac{{\partial \Xi }}{{\partial x}} = \dfrac{{\partial \Xi }}{{\partial X}}\dfrac{{\partial X}}{{\partial x}} + \dfrac{{\partial \Xi }}{{\partial Y}}\dfrac{{\partial Y}}{{\partial x}}$, and $\dfrac{{\partial \Xi }}{{\partial y}} = \dfrac{{\partial \Xi }}{{\partial X}}\dfrac{{\partial X}}{{\partial y}} + \dfrac{{\partial \Xi }}{{\partial Y}}\dfrac{{\partial Y}}{{\partial y}}$.\\

The components of the matrix $\mathbf{G}$ are given by

\begin{equation}
     {\bf{G}}_I^{std} = \left[ {\begin{array}{*{20}{c}}
     	\vspace{1.5mm}{\dfrac{{\partial {N_I}}}{{\partial x}}}&0\\
     	\vspace{1.5mm}0&{\dfrac{{\partial {N_I}}}{{\partial x}}}\\
     	\vspace{1.5mm}{\dfrac{{\partial {N_I}}}{{\partial y}}}&0\\
     	0&{\dfrac{{\partial {N_I}}}{{\partial y}}}
     	\end{array}} \right]{\rm{ }}
     \label{EqB4}
\end{equation}\\
\begin{equation}
{\bf{G}}_I^H = \left[ {\begin{array}{*{20}{c}}
	\vspace{1.5mm}{\dfrac{{\partial ({N_I}\bar H({\bf{x}}))}}{{\partial x}}}&0\\
	\vspace{1.5mm}0&{\dfrac{{\partial ({N_I}\bar H({\bf{x}}))}}{{\partial x}}}\\
	\vspace{1.5mm}{\dfrac{{\partial ({N_I}\bar H({\bf{x}}))}}{{\partial y}}}&0\\
	0&{\dfrac{{\partial ({N_I}\bar H({\bf{x}}))}}{{\partial y}}}
	\end{array}} \right]{\rm{ }}
\label{EqB5}
\end{equation}\\
\begin{equation}
{\bf{G}}_I^{Tip} = \left[ {\begin{array}{*{20}{c}}
	\vspace{1.5mm}{\dfrac{{\partial ({N_I}{{\bar A}_j}({\bf{X}})R({\bf{X}}))}}{{\partial x}}}&0\\
	\vspace{1.5mm}0&{\dfrac{{\partial ({N_I}{{\bar A}_j}({\bf{X}})R({\bf{X}}))}}{{\partial x}}}\\
	\vspace{1.5mm}{\dfrac{{\partial ({N_I}{{\bar A}_j}({\bf{X}})R({\bf{X}}))}}{{\partial y}}}&0\\
	0&{\dfrac{{\partial ({N_I}{{\bar A}_j}({\bf{X}})R({\bf{X}}))}}{{\partial y}}}
	\end{array}} \right]
\label{EqB6}
\end{equation}

%\section*{Acknowledgments}

%\section*{References}
%\bibliographystyle{model3-num-names}
%\bibliographystyle{elsarticle-num}
%\bibliographystyle{plain}
%\bibliographystyle{plainmat}
%\bibliography{HaiBib}

\begin{thebibliography}{10}
\expandafter\ifx\csname url\endcsname\relax
  \def\url#1{\texttt{#1}}\fi
\expandafter\ifx\csname urlprefix\endcsname\relax\def\urlprefix{URL }\fi
\expandafter\ifx\csname href\endcsname\relax
  \def\href#1#2{#2} \def\path#1{#1}\fi

\bibitem{Belytschko_Black_1999}
T.~Belytschko, T.~Black, Elastic crack growth in finite elements with minimal
  remeshing, International Journal for Numerical Methods in Engineering 45~(5)
  (1999) 601--620.

\bibitem{Moes_Dolbow_1999}
N.~Mo\"es, J.~Dolbow, T.~Belytschko, A finite element method for crack growth
  without remeshing, International Journal for Numerical Methods in Engineering
  46~(3) (1999) 131--150.

\bibitem{Daux_Moes_2000}
C.~Daux, N.~Mo\"es, J.~Dolbow, N.~Sukumar, T.~Belytschko, Arbitrary branched
  and intersecting cracks with the extended finite element method,
  International Journal for Numerical Methods in Engineering 48~(12) (2000)
  1741--1760.

\bibitem{Melenk_Babuska_1996}
J.~M. Melenk, I.~Babuska, The partition of unity finite element method: Basic
  theory and applications, Computer Methods in Applied Mechanics and
  Engineering 139~(1-4) (1996) 289--314.

\bibitem{Pant_Singh_Mishra_2010}
M.~Pant, I.~V. Singh, B.~K. Mishra, Numerical simulation of thermo-elastic
  fracture problems using element free {Galerkin} method, International Journal
  of Mechanical Sciences 52~(12) (2010) 1745--1755.

\bibitem{Ghorash_Valizadeh_Mohammad_2011i}
S.~S. Ghorashi, N.~Valizadeh, S.~Mohammadi, Extended isogeometric analysis for
  simulation of stationary and propagating cracks, International Journal for
  Numerical Methods in Engineering 89~(9) (2011) 1069--1101.

\bibitem{Duflot_2008}
M.~Duflot, The extended finite element method in thermoelastic fracture
  mechanics, International Journal for Numerical Methods in Engineering 74~(5)
  (2008) 827–847.

\bibitem{Fries_Belytschko_2010}
T.~P. Fries, T.~Belytschko, The extended/generalized finite element method: an
  overview of the method and its application, International Journal for
  Numerical Methods in Engineering 84~(3) (2010) 253--304.

\bibitem{Bayesteh_Mohammadi_2011}
H.~Bayesteh, S.~Mohammadi, {XFEM} fracture analysis of shell: The effect of
  crack tip enrichments, Computational Materials Science 50~(10) (2011)
  2793–2813.

\bibitem{Mariani_Perego_2003}
S.~Mariani, U.~Perego, Extended finite element method for quasi-brittle
  fracture, International Journal for Numerical Methods in Engineering 58~(11)
  (2003) 103–126.

\bibitem{Unger_Eckardt_Konke}
J.~Unger, S.~Eckardt, C.~Konke, Modelling of cohesive crack growth in concrete
  structures with the extended finite element method, Computer Methods in
  Applied Mechanics and Engineering 196~(41-44) (2007) 4087–4100.

\bibitem{Richardson_Hegemann}
C.~Richardson, J.~Hegemann, E.~Sifakis, J.~Hellrung, J.~Teran, An {XFEM} method
  for modeling geometrically elaborate crack propagation in brittle materials,
  International Journal for Numerical Methods in Engineering 88~(10) (2011)
  1042--1065.

\bibitem{Baydoun_Fries_2012}
M.~Baydoun, T.~Fries, Crack propagation criteria in three dimensions using the
  {XFEM} and an explicit – implicit crack description, International Journal
  of Fracture 179~(1-2) (2012) 51--70.

\bibitem{Pedro_Belytschko_2005}
M.~Pedro, T.~Belytschko, Analysis of three-dimensional crack initiation and
  propagation using the extended finite element method, International Journal
  for Numerical Methods in Engineering 63~(5) (2005) 760--788.

\bibitem{Sukumar_Tabarraei_2004}
N.~Sukumar, A.~Tabarraei, Conforming polygonal finite elements, International
  Journal for Numerical Methods in Engineering 61~(12) (2004) 2045–2066.

\bibitem{Sukumar_Malsch_2006}
N.~Sukumar, E.~A. Malsch, Recent advances in the construction of polygonal
  finite element interpolants, Archives of Computational Methods in Engineering
  13 (2006) 129–163.

\bibitem{Talischi_Paulino_Pereira_2012}
C.~Talischi, G.~H. Paulino, A.~Pereira, I.~F.~M. Menezes, Polymesher: a
  general-purpose mesh generator for polygonal elements written in matlab,
  Structural and Multidisciplinary Optimization 45~(3) (2012) 309–328.

\bibitem{Talischi_Pereira_Paulino_2014}
C.~Talisch, A.~Pereira, G.~H. Paulino, I.~F.~I. Menezes, M.~S. Carvalho,
  Polygonal finite elements for incompressible fluid flow, International
  Journal for Numerical Methods in Fluids 74~(2) (2014) 134–151.

\bibitem{Biabanaki_2014}
S.~O.~R. Biabanaki, A.~R. Khoei, P.~Wriggers, Polygonal finite element methods
  for contact-impact problems on non-conformal meshes, Computer Methods in
  Applied Mechanics and Engineering 269 (2014) 198--221.

\bibitem{Khoei_2015}
A.~R. Khoei, R.~Yasbolaghi, S.~O.~R. Biabanaki, A polygonal-{FEM} technique in
  modeling large sliding contact on non-conformal meshes: a study on polygonal
  shape functions, Engineering Computations 32~(5) (2015) 1391–1431.

\bibitem{Gain_Talischi_Paulino_2014}
A.~L. Gain, C.~Talischi, G.~H. Paulino, On the virtual element method for
  three-dimensional linear elasticity problems on arbitrary polyhedral meshes,
  Computer Methods in Applied Mechanics and Engineering.

\bibitem{Ghosh_Moorthy_1995}
S.~Ghosh, S.~Moorthy, Elastic-plastic analysis of arbitrary heterogeneous
  materials with the voronoi cell finite element method, Computer Methods in
  Applied Mechanics and Engineering 121~(1-4) (1995) 373–409.

\bibitem{Chi_Talischi_Pamies_Paulino_2014}
H.~Chi, C.~Talischi, O.~Lopez-Pamies, G.~H. Paulino, Polygonal finite elements
  for finite elasticity, International Journal for Numerical Methods in
  Engineering 101~(4) (2014) 305–328.

\bibitem{Ooi_2017}
E.~T. Ooi, C.~Song, S.~Natarajan, A scaled boundary finite element formulation
  with bubble functions for elasto-static analyses of functionally graded
  materials, Computational Mechanics 60~(6) (2017) 943--967.

\bibitem{Pramod_2018}
A.~L.~N. Pramoda, E.~T. Ooi, C.~Song, S.~Natarajan, Numerical estimation of
  stress intensity factors in cracked functionally graded piezoelectric
  materials - {A} scaled boundary finite element approach, Composite Structures
  206 (2018) 301--312.

\bibitem{Chi_2017}
H.~Chi, L.~B. da~Veiga, G.~H. Paulino, Some basic formulations of the {V}irtual
  {E}lement {M}ethod ({VEM}) for finite deformations, Computer Methods in
  Applied Mechanics and Engineering 318 (2017) 148--192.

\bibitem{Chi_2018}
H.~Chi, L.~B. da~Veiga, G.~H. Paulino, A simple and effective gradient recovery
  scheme and a posteriori error estimator for the {V}irtual {E}lement {M}ethod
  ({VEM}), Computer Methods in Applied Mechanics and Engineering\href
  {https://doi.org/https://doi.org/10.1016/j.cma.2018.08.014}
  {\path{doi:https://doi.org/10.1016/j.cma.2018.08.014}}.

\bibitem{Hung_2017}
H.~Nguyen-Xuan, S.~Nguyen-Hoang, T.~Rabczuk, K.~Hackl, A polytree-based
  adaptive approach to limit analysis of cracked structures, Computer Methods
  in Applied Mechanics and Engineering 313 (2017) 1006--1039.

\bibitem{Hung_2018}
H.~Nguyen-Xuan, H.~V. Do, K.~N. Chau, An adaptive strategy based on conforming
  quadtree meshes for kinematic limit analysis, Computer Methods in Applied
  Mechanics and Engineering 341 (2018) 485--516.

\bibitem{Nhon_Valizadeh_2015}
N.~Nguyen-Thanh, N.~Valizadeh, M.~N. Nguyen, H.~Nguyen-Xuan, et~al, An extended
  isogeometric thin shell analysis based on {Kirchhoff-Love} theory, Computer
  Methods in Applied Mechanics and Engineering 284 (2015) 265--291.

\bibitem{Ghorashi_Valizadeh_2015}
S.~S. Ghorashi, N.~Valizadeh, S.~Mohammadi, T.~Rabczuk, T-spline based {XIGA}
  for fracture analysis of orthotropic media, Computers {\&} Structures 147
  (2015) 138--146.

\bibitem{Zhuang_Augarde_Bordas_2011}
X.~Zhuang, C.~Augarde, S.~Bordas, Accurate fracture modelling using meshless
  methods, the visibility criterion and level sets: Formulation and {2D}
  modelling, International Journal for Numerical Methods in Engineering 86~(2)
  (2011) 249--268.

\bibitem{Tabarraei_Sukumar_2008}
A.~Tabarraei, N.~Sukumar, Extended finite element method on polygonal and
  quadtree meshes, Computer Methods in Applied Mechanics and Engineering
  197~(5) (2008) 425–438.

\bibitem{Khoei_Yasbo_2015}
A.~R. Khoei, R.~Yasbolaghi, S.~O.~R. Biabanaki, A polygonal finite element
  method for modeling crack propagation with minimum remeshing, International
  Journal of Fracture 194~(2) (2015) 123–148.

\bibitem{Philippe_Wolfgang_1994_1}
H.~G. Philippe, G.~K. Wolfgang, Finite strains at the tip of a crack in a sheet
  of hyperelastic material: {I}. {Homogeneous case}, Journal of Elasticity
  35~(1-3) (1994) 61--98.

\bibitem{Philippe_Wolfgang_1994_2}
H.~G. Philippe, G.~K. Wolfgang, Finite strains at the tip of a crack in a sheet
  of hyperelastic material: {II}. {Special bimaterial cases}, Journal of
  Elasticity 35~(1-3) (1994) 99--137.

\bibitem{Dolbow_Moes_Belytschko_2001}
J.~Dolbow, N.~Mo\"es, T.~Belytschko, An extended finite element method for
  modeling crack growth with frictional contact, Computer Methods in Applied
  Mechanics and Engineering 190~(51-52) (2001) 6825--6846.

\bibitem{Wells_Sluysi_2001}
G.~N. Wells, L.~J. Sluys, A new method for modelling cohesive cracks using
  finite elements, International Journal for Numerical Methods in Engineering
  50~(12) (2001) 2667–2682.

\bibitem{Rabczuk_Zi_2008}
T.~Rabczuk, G.~Zi, S.~Bordas, H.~Nguyen-Xuan, A geometrically non-linear
  three-dimensional cohesive crack method for reinforced concrete structures,
  Engineering Fracture Mechanics 75~(16) (2008) 4740--4758.

\bibitem{Dolbow_Devan_2004}
J.~E. Dolbow, A.~Devanl, Enrichment of enhanced assumed strain approximations
  for representing strong discontinuities: addressing volumetric
  incompressibility and the discontinuous patch test, International Journal for
  Numerical Methods in Engineering 59~(1) (2004) 47--67.

\bibitem{Legrain_Moes_Verron_2005}
G.~Legrain, N.~Mo\"es, E.~Verron, Stress analysis around crack tips in finite
  strain problems using the extended finite element method, International
  Journal for Numerical Methods in Engineering 63~(2) (2005) 290--314.

\bibitem{Karoui_2014}
A.~Karoui, K.~Mansouri, Y.~Renard, M.~Arfaoui, The extended finite element
  method for cracked hyperelastic materials: {A} convergence study,
  International Journal for Numerical Methods in Engineering 100~(3) (2014)
  222--242.

\bibitem{Rashetnia_2015}
R.~Rashetnia, S.~Mohammadi, Finite strain fracture analysis using the extended
  finite element method with new set of enrichment functions, International
  Journal for Numerical Methods in Engineering 102~(6) (2015) 1316--1351.

\bibitem{Broumand_Khoei_2013}
P.~Broumand, A.~R. Khoei, The extended finite element method for large
  deformation ductile fracture problems with a non-local damage-plasticity
  model, Engineering Fracture Mechanics 112-113 (2013) 97--125.

\bibitem{Wriggers_2008}
P.~Wriggers, Nonlinear Finite Element Methods, Springer, New York, NY, 2008.

\bibitem{Belytschko_Liu_Moran_2013}
T.~Belytschko, W.~Liu, B.~Moran, Nonlinear Finite Elements for Continua and
  Structures, Wiley, Chichester, West Sussex, UK, 2000.

\bibitem{Knowles_Sternberg_1973}
J.~K. Knowles, E.~Sternberg, An asymptotic finite-deformation analysis of the
  elastostatic field near the tip of a crack, Journal of Elasticity 3~(2)
  (1973) 67–107.

\bibitem{Fries_2008}
T.~P. Fries, A corrected {XFEM} approximation without problems in blending
  elements, International Journal for Numerical Methods in Engineering 75~(5)
  (2008) 503--532.

\bibitem{Wachspress_1979}
E.~L. Wachspress, Rational basis functions for curved elements, Computers and
  Mathematics with Applications 5~(4) (1979) 235–239.

\bibitem{Sibson_2012}
R.~Sibson, A vector identity for the {Dirichlet} tesselation, Mathematical
  Proceedings of the Cambridge Philosophical Society 87~(1) (1980) 151–155.

\bibitem{Floater_2014}
M.~S. Floater, Wachspress and mean value coordinates, Approximation Theory XIV:
  San Antonio 2013, Springer (2014) 81--102.

\bibitem{Floater_2016}
M.~S. Floater, Generalized barycentric coordinates and applications, Acta
  Numerica 24 (2015) 161–214.

\bibitem{Hormann_2006}
K.~Hormann, M.~S. Floater, Mean value coordinates for arbitrary planar
  polygons, Advances in Computational Mathematics 25~(4) (2006) 1424--1441.

\bibitem{Talischi_2015}
C.~Talischi, A.~Pereira, I.~F.~M. Menezes, G.~H. Paulino, Gradient correction
  for polygonal and polyhedral finite elements, International Journal for
  Numerical Methods in Engineering 102~(3-4) (2015) 728--747.

\bibitem{Chi_2016}
H.~Chi, C.~Talischib, O.~Lopez-Pamies, G.~H. Paulino, A paradigm for
  higher-order polygonal elements in finite elasticity using a gradient
  correction scheme, Computer Methods in Applied Mechanics and Engineering 306
  (2016) 216--251.

\bibitem{Gupta_1978}
Gupta, A finite element for transition from a fine grid to a coarse grid,
  International Journal for Numerical Methods in Engineering 12~(1) (1978)
  35--45.

\bibitem{Tabarraei_2007}
A.~Tabarraei, N.~Sukumar, Adaptive computations using material forces and
  residual-based error estimators on quadtree meshes, Computer Methods in
  Applied Mechanics and Engineering 196~(25-28) (2007) 2657–2680.

\bibitem{Lake_1970}
G.~J. Lake, Application of fracture mechanics to failure in rubber articles,
  with particulary reference to groove cracking in tyres, International
  Conference on Deformation, Yield and Fracture of Polymers, Cambridge.

\bibitem{Lindley_1972}
P.~B. Lindley, Energy for crack growth in model rubber components, The Journal
  of Strain Analysis for Engineering Design 7~(2) (1972) 132--140.

\bibitem{Yeoh_2002}
O.~H. Yeoh, Relation between crack surface displacements and strain energy
  release rate in thin rubber sheets, Mechanics of Materials 34~(8) (2002)
  459--474.

\bibitem{Steinmann_Ackermann_Barth_2001}
P.~Steinmann, D.~Ackermann, F.~J. Barth, Application of material forces to
  hyperelastostatic fracture mechanics. {II}. {Computational setting},
  Computational setting. International Journal of Solids and Structures
  38~(32-33) (2001) 5509--5526.

\end{thebibliography}
  
 \newcommand{\noop}[1]{}

\clearpage

\end{document}